\documentclass[11pt]{article}
\usepackage{amsmath}
\usepackage{amsthm}
\usepackage{amsfonts,mathrsfs}
\usepackage{amssymb,latexsym}
\usepackage{amsbsy}
\usepackage{amsxtra}
\usepackage{graphicx}

\newtheorem{thm}{Theorem}[section]
\newtheorem{prop}[thm]{Proposition}
\newtheorem{lem}[thm]{Lemma}
\newtheorem{rem}[thm]{Remark}

\def\para#1{\vskip .4\baselineskip\noindent{\bf #1}}

\def\para#1{\vskip 0.4\baselineskip\noindent{\bf #1}}

\newcommand{\thmref}[1]{Theorem~{\rm \ref{#1}}}
\newcommand{\lemref}[1]{Lemma~{\rm \ref{#1}}}

\newcommand{\propref}[1]{Proposition~{\rm \ref{#1}}}

\def\qed{$\qquad \Box$}

\def\cd{(\cdot)}

\def\rr{{\Bbb R}}

\def\nd{\noindent}
\renewcommand{\epsilon}{\varepsilon}

\newcommand{\EE}{\mathbb{E}}
\newcommand{\PP}{\mathbb{P}}

\makeatletter \@addtoreset{equation}{section}

\newcommand{\beq}[1]{\begin{equation} \label{#1}}
\newcommand{\eeq}{\end{equation}}
\newcommand{\bed}{\begin{displaymath}}
\newcommand{\eed}{\end{displaymath}}

\newcommand{\bea}{\bed\begin{array}{rl}}
\newcommand{\eea}{\end{array}\eed}
\newcommand{\ad}{&\!\!\!\disp}
\newcommand{\aad}{&\disp}
\newcommand{\barray}{\begin{array}{ll}}
\newcommand{\earray}{\end{array}}
\newcommand{\disp}{\displaystyle}
\newcommand{\al}{\alpha}

\begin{document}

\title{On Laws of Large Numbers for Systems with Mean-Field Interactions and Markovian Switching}
\author{Son L. Nguyen,\thanks{Department of Mathematics,
University of Puerto Rico, Rio Piedras campus, San Juan, PR 00936, USA, sonluu.nguyen@upr.edu. This
research was supported by a seed fund of Department of Mathematics at
University of Puerto Rico, Rio Piedras campus.} \and George Yin,\thanks{Department of Mathematics, Wayne State University,
Detroit, MI 48202, USA, gyin@math.wayne.edu. This research was supported in part by the Army Research Office.} \and Tuan A. Hoang\thanks{Department of Mathematics, Wayne State University,
Detroit, MI 48202, USA, tuan.hoang@wayne.edu.}}

\maketitle

\begin{abstract}
Focusing on stochastic systems arising in mean-field models, the systems under consideration belong to the class of switching diffusions, in which continuous dynamics and discrete events coexist and interact. The discrete events are modeled by a continuous-time Markov chain. Different from the usual switching diffusions, the systems include mean-field interactions.
Our effort is devoted to obtaining laws of large numbers for the underlying systems. One of the distinct features of the paper is the limit of the empirical measures is not deterministic  but a random measure depending on the history of the Markovian switching process. A main difficulty is that the standard martingale approach cannot be used to characterize the limit because of the coupling due to the random switching process.  In this paper, in  contrast to the classical approach, the limit is characterized as the conditional distribution (given the history of the switching process) of the solution to a stochastic McKean-Vlasov  differential equation with  Markovian switching.

\medskip
\nd{\bf Key Words.} Mean-field model, Markovian switching process, law of large number, McKean-Vlasov equation.

\medskip
\nd{\bf Mathematics Subject Classification.} 60J25, 60J27, 60J60, 93E20.

\medskip
\nd{\bf Running title.} LLN for Systems with Mean-Field Interactions and Markovian Switching

\end{abstract}

\newpage

\section{Introduction}\label{sec:int}

This work
focuses on laws of large numbers for a class of stochastic systems involving mean-field interactions and random switching.  It
 is motivated by two lines of recent advances in the study of stochastic systems and applications. One of them is the emerging interests in mean-field models, and the other is the use of regime switching in stochastic systems.

Originated from statistical physics, mean-field models describe stochastic systems containing
a large number of particles having weak interactions. To overcome the complexity of interactions due to
the large scale of system, all interactions with each particle are replaced by a single averaged interaction
(naturally represented by an empirical measure associated to system).  One of the first mathematical treatments was the influential work of Dawson \cite{Dawson83}, which rigorously justified the replacement of  a large number of ``bodies'' by  a representative  ``body'' in many-body problems. Also obtained in that paper was phase transition properties.  The subsequent work \cite{DZ1991} delineates some limit theory for jump mean-field models. Although originally appeared
 in physics, mean-field models have arisen in many different application areas, including communication networks, mathematical finance, chemical and biological systems, and social sciences. For an extensive list of references to such applications, see \cite{BaladronFFT12,BFY13,ContucciGM08,Sznitman91}.
 Recently,
  renewed interest has been shown using
 mean-field models in game theory, which is originated independently in the work of  Huang,  Caines, and Malham\'e \cite{HCM03,HMC06}, and
 Lasry and Lions \cite{LL06a}.  The mean-field interaction has
been used to model the weak interaction between players in large population games and the limiting
results are used to construct computable decentralized strategies, leading to substantial progress in the development of mean-field game theory.  The book by
 Bensoussan, Frehse, and Yam
 \cite{BFY13} provided an illuminating
  presentation and discussion of certain aspects of mean-field games and  mean-field type controls,
  describing their similarities and differences together with a unified approach for treating them.
  A more analytic approach can be found in
  the book by Kolokoltsov
  \cite{Kolokoltsov10}.
  It has been seen that
 mean-field models and mean-field games enjoy a wide range of
  applications and potential applications in
  economics, social networks,
cyber physical systems, and other branches of sciences
and engineering; see \cite{BFY13,HN16,HCM03,HMC06,NH12,
WangZ12,WangZ17} and references therein.

Along another line, the so-called hybrid systems have gained increasing popularity due to their ability to handle numerous real-world applications in which discrete and continuous dynamics coexist and interact. One class of such hybrid systems is switching diffusions. Take for instance, applications in control systems and optimization.
One of the most widely used control engineering models in the literature
is the so-called linear quadratic Gaussian regulator problem;
see \cite{FlemingR} for a
traditional model. For many new applications in networked systems, it has been found that in addition to the random noise represented by Brownian type of disturbances, there is a source of randomness owing to the presence of random environment
that displays pure jump behavior and that can be modeled by a continuous-time Markov chain. As a result, one has a controlled switching diffusion instead of controlled diffusion as in the traditional setup; see \cite{YinZhu09} for some recent results on switching diffusions and applications.
For a wide variety of applications, we mention
the work on flexible manufacturing systems \cite{SethiZ},
approximation to invariant measures \cite{BSY16},
controlled piecewise deterministic Markov processes \cite{CD13},
population dynamics \cite{LM09},
business cycle models in random environment \cite{TYW16},
stochastic approximation \cite{YKI04} with applications to wireless communication such as spreading code optimization and adaptation in CDMA,
Markowitz mean-variance portfolio selection \cite{ZhY03}, and
Lotka-Volterra models in ecology \cite{ZY09}, among others.
Furthermore, there has been new effort in treating switching diffusions in conjunction with mean-field interactions \cite{XiYin09}.

Motivated by the aforementioned two aspects,
the focus of the current paper lies in the intersection of the mean-field models and the switching diffusion models. We concentrate on large-scale systems with weak interactions in a random environment represented by switching diffusions
in which the  Markov chains delineate random environment
changes not represented by the usual diffusions.
Recently, some related works  have been considered in
\cite{WangZ12,WangZ17} for studying
mean-field games and social optimality. In this work, we investigate functional laws of large numbers for such systems.

Why should we be so concerned about laws of large numbers and why should such an effort be necessary? Not only is the study interesting from a mathematical point of view, but also it is crucial from a practical consideration. Treating large-scale systems, a main effort is to reduce the computation complexity. Laws of large numbers provide us with an effective machinery to overcome the difficulties.
As a motivational example, consider a mean-field game problem with $N$ players for a large number $N$.
Let $x_i (t) \in \rr^d$, $1\le i\le N$, be the state of player $i$ that satisfies the following equation
\bea dx_i(t)\ad = b\bigg(x_i(t),{1\over N}\sum_{j=1}^N\delta_{x_j(t)},\alpha(t_-), u_i(t)\bigg)dt + \sigma\bigg(x_i(t),{1\over N}\sum_{j=1}^N\delta_{x_j(t)},\alpha(t_-)\bigg) dw_i(t)\\
   x_i(0)\ad = x_i,  \ i\le N, \ \al(0)=\al.\eea
where
$b(\cdot,\cdot,\cdot,\cdot)$ and $\sigma(\cdot,\cdot,\cdot)$ are appropriate functions,
$w_i\cd$, $1\le i\le N$, are independent $\rr^d$-valued Brownian motions,
$\al(\cdot)$ is  a continuous-time Markov chain independent of the Brownian motions $w_i(\cdot)$, $\delta_x(\cdot)$ denotes the Dirac measure centered at $x$ for each $x \in \rr^d$, and $u_i\cd$, $1\le i\le N$, is the control of player $i$ taking values in a compact subset of another Euclidean space $\rr^{d_1}$.  Player $i$, $1\le i\le N$, wishes to
minimize its own cost
$$J_i(x_j,u_j(\cdot): j\le N)
= \EE_{\{x_j: j\le N\}}\int^T_0 R\bigg(x_i(t),{1\over N}\sum_{j=1}^N\delta_{x_j(t)}, u_i(t)\bigg) dt ,$$
where  $R(\cdot,\cdot,\cdot)$ is a running cost function and the expectation is taken with $x_j(0)=x_j$.
To obtain low complexity strategies, consistent mean-field approximations provide a powerful approach. Consequently,
each player only needs to know its own state information and the aggregate effect of the overall population, which may be pre-computed off-line. A crucial step of this approach is to approximate the instantaneous measure ${1\over N}\sum_{j=1}^N\delta_{x_j(t)}$ of the processes under consideration by a stationary measure as $N\to \infty$. In order to take such a step, one needs to demonstrate that the system indeed possesses such a limit measure.   The law of large numbers of the corresponding systems provides the existence of this limit and helps to characterize it. With the motivation for finding optimal strategies for
 mean-field models with $N$ players and Markovian switching, this work establishes the laws of large numbers for such systems and paves a way for solving the underlying problem.

Regarding law of large numbers, it is worth mentioning that since the pioneering works of Kac \cite{Kac56} and McKean \cite{McKean66}, many important results have been obtained for investigating the time evolution of stochastic systems with long range weak interactions.
Many variants of such systems have also been examined.
For example, in \cite{KurtzXiong99}, limit theorems were
established for a model in which there is a common space noise process that influences the dynamics of
each particle. Law of large numbers  in a setting where particle evolution
depends on independent jumps and switching processes were studied in \cite{LDF17,Oelschlager84} and \cite{Graham90}, respectively. In \cite{DawsonV95}, law of large numbers was studied for a model where the noises are correlated.

One of the novel features of this paper is the limit of the empirical measures is not deterministic  but a random measure that depends on the history of the Markovian switching process. In addition, the stochastic McKean-Vlasov equation in the limit is driven by martingales associate with the Markov switching process.
 As a consequence, there is a main difficulty to characterize the limit using the
 martingale problem formulation as in \cite{DawsonV95,Gartner88,Graham90}. To overcome this difficulty,
we use a new approach.  Different from the classical work, we characterize
 the limit as the unique solution to a stochastic McKean-Vlasov equation with Markovian switching, which is represented by the conditional distribution of the solution to a McKean-Vlasov stochastic differential equation with a Markovian switching given the history of the switching process. In contrast,
 for  the problem treated in \cite{Graham90}, each particle possess its own switching process and the limit is represented as the distribution of solution to a
   McKean-Vlasov stochastic differential equation.
We note that in \cite{KurtzXiong99},
 Kurtz
and Xiong treated interacting particles.
In their paper, there is a common
space-time Gaussian white noise.  They obtained
 law of large numbers with the conditional distribution in the limit.
In their case, the martingale problem approach cannot be used either. Nevertheless, their model contains infinitely many exchangeable
particles. Thus ergodic theory can be applied
to the system, whereas in our case, we no longer have infinitely many exchangeable particles thus we cannot carry out the study by directly applying the existing ergodic theory.

In networked systems, the discrete component, namely, the random switching process often has a rather large state space. The transition among the states are not of the same speed. Some of them vary rapidly, whereas the others evolve slowly. As illustrated in \cite{SethiZ} (see also \cite{YinZhang13}), the applications demand the consideration of the so-called nearly decomposable structures.
Here nearly decomposable  is understood in the sense that the switching among different subspaces are still possible although they appear relatively infrequently.
 Consequently, the large state space is naturally divisible into a number of subspaces so that the transitions in each subspace take place
 at a fast pace; the transitions from one subspace to another
occur slowly. Such a situation leads to the modeling using
two-time scales
 as in \cite{YinZhang13} by introducing a small parameter $\epsilon>0$ into the systems.
 In this paper,
we will also investigate this case. The goal is still to get laws of large numbers. However, in lieu of one parameter $N$, we have two parameters $N$ and $\epsilon$. The limit is taken to be as $\epsilon\to 0$, $N\to \infty$, and $(1/\epsilon)\wedge N\to \infty$.

The rest of the paper is arranged as follows. Section \ref{sec:for} presents the formulation of the problem that we wish to study.
Section \ref{sec:pre} collects a number of preliminary results of interacting particle systems with Markovian switching.  Section \ref{sec:lln} demonstrates the law of large numbers for the systems. Section \ref{sec:two-time} examines systems in which the random switching displays two-time-scale behavior.
Finally, an appendix containing the proofs of some technical lemmas is placed at the end of the paper.
We remark that this paper is devoted to convergence in the form of law of large numbers. The rates of convergence is an interesting topic for future research. In the literature, some of such attempts can be found in
\cite{Kolokoltsov10} using an analytic approach and  \cite{HMC06} using martingale-type estimates. For problems under the setting of this paper, because of the conditional distributions
usage, careful thoughts and considerations are needed to treat the rate of convergence issue.

\section{Formulation}\label{sec:for}
We consider a mean-field system of $N$ particles (with $N$ being a large number),
described by the following system of stochastic differential equations
\begin{equation}\label{def-dynamics}
    dx_i(t) = b\bigg(x_i(t),{1\over N}\sum_{j=1}^N\delta_{x_j(t)},\alpha(t_-)\bigg)dt + \sigma\bigg(x_i(t),{1\over N}\sum_{j=1}^N\delta_{x_j(t)},\alpha(t_-)\bigg) dw_i(t),
\end{equation}
for $i=1,2,\ldots,N$, where $\delta_x(\cdot)$ denotes the Dirac measure centered at $x$ with $x\in\mathbb{R}^d$, $w_1(\cdot)$, $w_2(\cdot)$, $\ldots$, $w_N(\cdot)$ are $N$ independent $d$-dimensional standard Brownian motions, and $\alpha(\cdot)$ is a Markov chain taking values in a finite state space $\mathbb{S} =
\{1, 2, \ldots, m_0\}$ with a generator $Q=\big(q_{i_0j_0}\big)_{i_0,j_0\in\mathbb{S}}$ satisfying the following properties: $q_{i_0j_0}\ge0$ for $i_0\ne j_0\in\mathbb{S}$ and $q_{i_0i_0}=-\sum_{j_0\ne i_0}q_{i_0j_0}$ for each $i_0 \in {\mathbb S}$.

Throughout this paper, we assume that the Brownian motions $w_i(\cdot)$, $1\le i\le N$, and the Markov chain $\alpha(\cdot)$ are independent and defined on a common complete probability space $(\Omega,\mathcal{F}, \PP)$.
Note that the transition rule of the Markov chain $\alpha(t)$ satisfies
$$
\PP\Big( \alpha(t + \Delta t) = j_0 \Big| \alpha (t) = i_0, \
0\leq s \leq t \Big) = q_{i_0j_0}\Delta t + o(\Delta t),
$$
for any pair $i_0,j_0\in\mathbb{S}$. It is clear that $x_i(t)$ depends on $N$  in accordance with \eqref{def-dynamics}, but
 to simplify the notation,  we omit the index $N$ in $x_i(t)$ in what follows.

\para{Notation.}
Let $C_b(\mathbb{R}^d)$ denote the space of bounded and continuous functions on $\mathbb{R}^d$ equipped with the usual supremum norm $\|\cdot\|$, $C^k_b(\mathbb{R}^d)$ the space of all functions in $C_b(\mathbb{R}^d)$ whose partial derivatives up to order $k$
 are bounded and continuous,
and  $C^k_c(\mathbb{R}^d)$ the space of functions whose partial
  derivatives up to order $k$ are continuous with compact support.
Denote by $\mathscr{M}_1$ the space of all probability measures on $\mathbb{R}^d$. For $f\in C_b(\mathbb{R}^d)$ and $\mu\in\mathscr{M}_1$, define $\langle \mu,f\rangle=\int_{\mathbb{R}^d}f(x)\mu(dx)$.   We shall use the total variation metric $\|\cdot\|_{TV}$ and the bounded Lipschitz metric  $\|\cdot\|_{BL}$ on $\mathscr{M}_1$ given as
$$
\|\mu-\eta\|_{BL}=\sup\Bigg\{\big|\big\langle\mu,f\big\rangle-\big\langle\eta,f\big\rangle\big|: \|f\|\le1,\sup_{x\ne y\in\mathbb{R}^d}{|f(x)-f(y)|\over|x-y|}\le1\Bigg\},
$$
for $\mu,\eta\in\mathscr{M}_1$.
 It follows from \cite{Dudley66} that $(\mathscr{M}_1,\|\cdot\|_{BL})$ is a separable and complete metric space, which is topologically equivalent to the space of all probability measures on $\mathbb{R}^d$ equipped with the weak topology. Endow $\mathbb{S}$ with a metric $d_\mathbb{S}$ satisfying $d_\mathbb{S}(i_0,i_0)=0$ and $d_\mathbb{S}(i_0,j_0)=1$ if $i_0\ne j_0$ for $i_0,j_0\in\mathbb{S}$.  Define the following metric $d$ on the product space $\mathscr{M}_1\times\mathbb{S}$,
\begin{equation}\label{Def-metric-d}
d\big((\mu,i_0),(\eta,j_0)\big)= \big\|\mu-\eta\big\|_{BL}+d_\mathbb{S}\big(i_0,j_0\big),\quad \forall\,\mu,\eta\in\mathscr{M}_1,i_0,j_0\in\mathbb{S}.
\end{equation}

For a metric space $E$, let $\mathcal{B}(E)$ be the Borel $\sigma$-field on $E$ and $\mathcal{P}(E)$ denote the space of all probability measures on $\big(E,\mathcal{B}(E)\big)$ equipped with the weak convergence topology.  Let $C([0,T],E)$ denote the space of all continuous functions $h: [0,T]\to E$ equipped with the supremum metric and $D([0,T],E)$  the space of all c\`adl\`ag
functions $h: [0,T]\to E$ equipped with the usual Skorohod topology.  Denote by $D_f([0,T],\mathbb{S})$ the subspace of $D([0,T],\mathbb{S})$ which contains all processes with finite jumps. Since $\mathbb{S}$ is a discrete set, $D_f([0,T],\mathbb{S})$ is a closed subset of $D([0,T],\mathbb{S})$. For a given $\mu \in  \mathscr{M}_1$ and functions $f(\cdot,\cdot,\cdot)$ and $g(\cdot,\cdot)$ satisfying $f(\cdot,\cdot,i_0)\in C_b(\mathbb{R}\times\mathbb{R}^d)$  and $g(\cdot,i_0)\in C_b(\mathbb{R}^d)$ for each $i_0\in\mathbb{S}$,  we define $\langle\mu ,f(t,\cdot,i_0)\rangle=\int_{\mathbb{R}^d}f(t,x,i_0)\mu(dx)$ and $\langle\mu,g(\cdot,i_0)\rangle=\int_{\mathbb{R}^d}g(x,i_0)\mu(dx)$.

Let $\mathcal{B}(\mathbb{R}^d)$ denote the usual Borel $\sigma$-field on $\mathbb{R}^d$. For any vector $x\in\mathbb{R}^d$ or matrix $A\in\mathbb{R}^{d\times d}$, $|x|$ and $|A|$  denote their usual norms in $\mathbb{R}^d$ and $\mathbb{R}^{d\times d}$, respectively,  and $x'$ and $A'$  denote their transposes. In addition, the inner product of two vectors $x,y$ is denoted by $(x,y)$. In what follows, we frequently use two particular functions $\varphi(\cdot),\psi(\cdot):\mathbb{R}^d\to\mathbb{R}$ defined by $\varphi(x)=|x|$ and $\psi(x)=|x|^2$, respectively.
For $t>0$, denote $\mathcal{F}^\alpha_{t_-}=\sigma\big\{\alpha(s):0\le s<t\big\}$ and
$$
\mathcal{F}^{N,\alpha}_{t}=\sigma\big\{w_i(s),\alpha(s):0\le s\le t,1\le i\le N\big\}.
$$

For a random variable $\varsigma$ on $\big(\Omega,\mathcal{F},\mathbb{P}\big)$, we denote by \begin{equation}\label{eq:L-law}
\mathscr{L}(\varsigma) \hbox{ its  distribution and } \eta_t=\mathscr{L}\big(\varsigma\big|\mathcal{F}^\alpha_{t_-}\big) \hbox{ its conditional distribution given }\mathcal{F}^\alpha_{t_-}\end{equation} in the sense that
\begin{equation}\label{eq:L-law-1}
\mathbb{E}\big(f(\varsigma)\big|\mathcal{F}^\alpha_{t_-}\big)=
\int_{\mathbb{R}^d}f(x)\eta_t(dx) \hbox{ for any } f\in C_b(\mathbb{R}^d).\end{equation}
We make the following assumptions.

\para{Assumption A.}
\begin{itemize}
\item[(A1)] For each $i_0\in\mathbb{S}$, $b(\cdot,\cdot,i_0):\mathbb{R}^d\times\mathscr M_1\to\mathbb{R}^d$ and $\sigma(\cdot,\cdot,i_0):\mathbb{R}^d\times\mathscr M_1\to\mathbb{R}^{d\times d}$ are Lipschitz continuous in that, there is a constant $L$ such that
    $$
    \Big|b\big(x,\mu,i_0\big)-b\big(y,\eta,i_0\big)\Big|+\Big|\sigma\big(x,\mu,i_0\big)-\sigma\big(y,\eta,i_0\big)\Big|\le L\Big(\big|x-y\big|+\big\|\mu-\eta\big\|_{BL}\Big),
    $$
    for all $x,y\in\mathbb{R}^d$ and $\mu,\eta\in \mathscr M_1$.

\item[(A2)]  The $\mathbb{R}^d$-valued function $b(\cdot,\cdot,\cdot)$ satisfies
 $$
 \Big|b\big(x,\mu,i_0\big)\Big|\le C\Big(1+\big|x\big|+\big\langle\mu,\varphi\big\rangle\Big),\quad  (x,\mu,i_0)\in\mathbb{R}^d\times\mathscr{M}_1\times\mathbb{S},
 $$
for some constant $C$ and $\varphi:\mathbb{R}^d\to\mathbb{R}$, $\varphi(x)=|x|$ and the matrix-valued function $\sigma(\cdot,\cdot,\cdot)$ is bounded.
\end{itemize}

Note that in the above and throughout the paper, for notational simplicity, the same notion $|\cdot|$ is used to denote different norms in $\mathbb{R}^d$, or $(\mathbb{R}^d)^N$, or $\mathbb{R}^{d\times d}$. It should, however, be clear from the context which norm is currently used.

It will be shown in the next section that under Assumption (A), for each fixed $N\ge1$,  system \eqref{def-dynamics} has a unique solution $\big(x_1(t),x_2(t),\ldots,x_N(t)\big)$. For $N\ge1$, $0\le t\le T$, and $A\in\mathcal{B}(\mathbb{R}^d)$, define
\begin{equation}\label{def-mu_N}
\mu_{N}(t,A) = \frac{1}{N}\sum_{j = 1}^N \delta_{x_j(t)}(A).
\end{equation}
Then $\mu_{N}(t,\cdot)$ is
a measured-valued process,
taking value on the space  $\mathscr{M}_1$ of probability measures on $\mathbb{R}^d$.
We denote by $\mathscr P_N$  the induced probability measure of $\big(\mu_N(\cdot),\alpha(\cdot)\big)$ on $D\big([0,T],\mathscr{M}_1\times\mathbb{S}\big)$. It can be shown that $\mathscr P_N$
concentrates on the set $C\big([0,T],\mathscr{M}_1\big)\times D_f\big([0,T],\mathbb{S}\big)$, a closed subspace of $D\big([0,T],\mathscr M_1\times\mathbb{S}\big)$.
  Using the notation mentioned thus far, in particular,  \eqref{eq:L-law} and \eqref{eq:L-law}, we proceed to derive the following main result. The proof is provided in Section \ref{sec:lln}, and some preliminary results are given in the next section as preparation.

\begin{thm}\label{thm-main1}
Assume {\rm(A1)}, {\rm(A2)}, and
$$
\sup_{N\in\mathbb N}\EE\big\langle \mu_N(0),\psi\big\rangle<\infty,\quad \mathscr{L}\big(\mu_N(0)\big)\Rightarrow \delta_{\mu_0}\text{ in } \mathcal{P}\big(\mathscr{M}_1,\|\cdot\|_{BL}\big),
$$
where $\psi:\mathbb{R}^d\to\mathbb{R}$ with $\psi(x)=|x|^2$. Then $\big(\mu_N(\cdot),\alpha(\cdot)\big)$ converges weakly to a process $\big(\mu_\alpha(\cdot),\alpha(\cdot)\big)$, where
$$
\big(\mu_\alpha(t),\alpha(t)\big)=\big(\mathscr L\big(y(t)\big|\mathcal{F}^\alpha_{t_-}\big),\alpha(t)\big),\quad 0\le t\le T,
$$
and $y(t)$, $0\le t\le T$, is the unique solution  of the following stochastic differential equation
$$
dy(t)=b\Big(y(t),\mathscr L\big(y(t)\big|\mathcal{F}^\alpha_{t_-}\big),\alpha(t_-)\Big)dt+\sigma\Big(y(t),\mathscr L\big(y(t)\big|\mathcal{F}^\alpha_{t_-}\big),\alpha(t_-)\Big)d\tilde w(t),\quad \mathscr L\big(y(0)\big)=\mu_0,
$$
where $\tilde w(\cdot)$ is a standard Brownian motion independent of $\alpha(\cdot)$.
\end{thm}

As mentioned in the introduction, motivated by applications in networked systems where the random switching process often has a large state space and the transition among the states are not at the same speed, we also treat mean-field systems that capture
different transition rates (slow and fast) of the switching process by using two-time scale approach. A parameter $\epsilon$ will be used to depict the difference of transition speeds. It can be shown that the law of large numbers also holds true for this case under some mild conditions similar to those in \thmref{thm-main1}. For clarity of presentation, the formulation of this case will be given in Section \ref{sec:two-time}.


\section{Preliminaries}\label{sec:pre}
In this section, we provide some preliminary results on weakly interacting systems with Markovian switching. For convenience, we first consider the general switching systems consisting of $N$-particles in $\rr^d$ without the weak interaction assumption. These systems can be formulated as switching diffusion processes in the larger space $(\rr^d)^N$. Weakly interacting systems of $N$-particles is then presented in Section 3.2 as a special case.

\subsection{General $N$-Particle System with Markovian Switching}

Let $x_{0,i}, 0\le i\le N$ be $\mathbb{R}^d$-valued random variables defined on $(\Omega,\mathcal{F},P)$ that are independent of $w_i(\cdot), 1\le i\le N,$ and $\alpha(\cdot)$. Assume that $\mathbb{E}|x_{0,i}|^2<\infty$ for $1\le i\le N$.  Consider the following stochastic differential equations with Markovian switching
\begin{align}
dx_i(t)& =\underline{b}_i\big(t,\underline{x}(t),\alpha(t_-)\big)dt+\underline{\sigma}_{i}\big(t,\underline{x}(t),\alpha(t_-)\big)dw_i(t),\quad  1\le i\le N,0\le t\le T,\label{Eq-SDE-MSwit-a}\\
x_i(0)& =x_{0,i} \quad\text{a.s.,}\notag
\end{align}
where $\underline{x}(t)=\underline{x}_N(t)=\big(x_1(t),x_2(t),\ldots,x_N(t)\big)\in\big(\mathbb{R}^d\big)^N$, $\underline{b}_i(\cdot,\cdot,\cdot):\mathbb{R}\times\big(\mathbb{R}^d\big)^N\times\mathbb{S}\to\mathbb{R}^d$, $\underline{\sigma}_{i}(\cdot,\cdot,\cdot):\mathbb{R}\times\big(\mathbb{R}^d\big)^N\times\mathbb{S}\to\mathbb{R}^{d\times d}$ are vector-valued functions, $w_1(\cdot),w_2(\cdot),\ldots,w_N(\cdot)$ are $\mathbb{R}^d$-valued independent standard Brownian motions, and $\alpha(\cdot)$ is a Markov chain with the state space $\mathbb{S}$ and generator $Q=\big(q_{i_0,j_0}\big)_{i_0,j_0\in\mathbb{S}}$ given as in the previous section.
Assume that for  $1\le i\le N$, $i_0\in\mathbb{S}$ and $0\le t\le T$, $\underline{b}_i\big(t,\cdot,i_0\big)$ and $\underline{\sigma}_{i}\big(t,\cdot,i_0\big)$ satisfy the following Lipschitz and linear growth conditions
\begin{equation}
\big|\underline{b}_i\big(t,\underline{x},i_0\big)-\underline{b}_i\big(t,\underline{y},i_0\big)\big|+\big|\underline{\sigma}_{i}\big(t,\underline{x},i_0\big)-\underline{\sigma}_{i}\big(t,\underline{y},i_0\big)\big|\le K\big|\underline{x}-\underline{y}\big|,\label{Cond-SDE-Lips}\\
\end{equation}
\begin{equation}
\big|\underline{b}_i\big(t,\underline{x},i_0\big)\big|+\big|\underline{\sigma}_{i}\big(t,\underline{x},i_0\big)\big|\le K\big(1+\big|\underline{x}\big|\big),\label{Cond-SDE-Linear}
\end{equation}
for any $\underline{x},\underline{y}\in (\mathbb{R}^d)^N$, where $K$  is a positive constant. It follows from Theorem 3.3.13 \cite{MaoYuan06} that under \eqref{Cond-SDE-Lips} and \eqref{Cond-SDE-Linear}, the system \eqref{Eq-SDE-MSwit-a} has a unique solution.

By virtue of \cite{MaoYuan06,YinZhu09}, for a function $V(\cdot,\cdot,\cdot):[0,T]\times(\mathbb{R}^d)^N\times\mathbb{S}\to\mathbb{R}$ such that  for each $i_0\in\mathcal{M}$, $V(\cdot,\cdot,i_0)\in C^{1,2}\big([0,T]\times(\mathbb{R}^d)^N\big)$, the generator of the general system of  $N$ particles is defined by
\begin{align}
\mathcal{L}_NV\big(t,\underline{x},i_0\big)&=\sum_{i=1}^N\underline{b}_i'\big(t,\underline{x},i_0\big)\nabla_{x_i}V\big(t,\underline{x},i_0\big)+{1\over2}\sum_{i=1}^N \Big(\underline{a}_{i}\big(t,\underline{x},i_0\big)\nabla_{x_i}\Big)'\nabla_{x_i}V\big(t,\underline{x},i_0\big)\notag\\
&\quad+\sum_{j_0\in\mathcal{M}}q_{i_0j_0}\Big(V\big(t,\underline{x},j_0\big)-V\big(t,\underline{x},i_0\big)\Big),\label{Def-Op-SDE-MSw}
\end{align}
for $\big(t,\underline{x},i_0\big)\in[0,T]\times(\mathbb{R}^d)^N\times\mathbb{S}$, where
$\nabla_{x_i}$ denotes the gradient with respect to $x_i$, and $\underline{a}_{i}\big(t,\underline{x},i_0\big)=\underline{\sigma}_{i}\big(t,\underline{x},i_0\big)\underline{\sigma}'_{i}\big(t,\underline{x},i_0\big)\in\mathbb{R}^{d\times d}$ for each $1\le i\le N$.

Associated with each pair $(i_0,j_0)\in {\mathbb S}\times {\mathbb S}$, $i_0\ne j_0$, the states of the Markov chain $\alpha(\cdot)$, define
\begin{equation}\label{Def-Quad-Var}
\big[M_{i_0j_0}\big](t)=\sum_{0\le s\le t}{1\!\!1}\big(\alpha(s_-)=i_0\big){1\!\!1}\big(\alpha(s)=j_0\big),\quad\quad \big\langle M_{i_0j_0}\big\rangle(t) =\int_0^tq_{i_0j_0}{1\!\!1}\big(\alpha(s_-)=i_0\big)ds,
\end{equation}
where ${1\!\!1}$ denotes the usual zero-one indicator function. It follows from Lemma IV.21.12 \cite{RogersWilliams00} that the process $M_{i_0j_0}(t)$, $0\le t\le T$, defined by
\begin{equation}\label{Def-Mart of MC}
M_{i_0j_0}(t)=\big[M_{i_0j_0}\big](t)- \big\langle M_{i_0j_0}\big\rangle(t)
\end{equation}
is a purely discontinuous and square integrable martingale with respect to $\mathcal{F}^{N,\alpha}_t$, which is null at the origin. The processes  $[M_{i_0j_0}](t)$ and $ \langle M_{i_0j_0}\rangle(t)$ are respectively its optional and predictable quadratic variations. For convenience, we define $M_{i_0i_0}(t)=\big[M_{i_0i_0}\big](t)=\big\langle M_{i_0i_0}\big\rangle(t)=0$ for each $i_0\in\mathbb{S}$.  From the definition of optional quadratic covariations (see Section 1.8 in \cite{LiptserShiryayev89}) we have the following orthogonality relation :
\begin{equation}\label{Eq-Ortho}
\big[w_i,w_j\big]=0\text{ when }i\ne j,\quad \big[M_{i_0j_0},w_j\big]=0,\quad\big[M_{i_0j_0},M_{p_0q_0}\big]=0 \text{ when }(i_0,j_0)\ne(p_0,q_0).
\end{equation}

For any function $V(\cdot,\cdot,\cdot):[0,T]\times(\mathbb{R}^d)^N\times\mathbb{S}\to\mathbb{R}$ such that $V(\cdot,\cdot,i_0)\in C^{1,2}\big([0,T]\times(\mathbb{R}^d)^N\big)$ for each $i_0\in\mathcal{M}$, we have the following It\^o formula
\begin{equation}\label{Eq-Ito-SDEMSw}
\begin{array}{ll}
V\big(t,\underline{x}(t),\alpha(t)\big)
&\!\! \disp =V\big(0,\underline{x}(0),\alpha(0)\big)+\int_0^t\bigg({\partial  \over\partial s}+\mathcal{L}_N\bigg)V\big(s,\underline{x}(s),\alpha(s_-)\big)ds\\
&\disp\quad+\sum_{i=1}^N\int_0^t\Big\langle\nabla_{x_i}V
\big(s,\underline{x}(s),\alpha(s_-)\big),
\underline{\sigma}_{i}\big(s,\underline{x}(s),
\alpha(s_-)\big)dw_i(s)\Big\rangle\\
&\disp \quad+\sum_{i_0\ne j_0}\int_0^t\Big(V\big(s,\underline{x}(s),
j_0\big)-V\big(s,\underline{x}(s),i_0
\big)\Big)dM_{i_0j_0}(s).
\end{array}\end{equation}
It can be seen that we have two martingales. One of them is due to the Brownian motion, whereas the other is
resulted from the jump process.

\subsection{$N$-Particle Mean-Field Model with Markovian Switching}
For $\underline{x}=\big(x_1,x_2,\ldots,x_N\big)\in(\mathbb{R}^d)^N$, denote the associated empirical probability measure $\delta_{\underline{x}}$ on $\mathscr{M}_1$ by
$\delta_{\underline{x}}={1\over N}\sum_{j=1}^N\delta_{x_j}$.
Consider the system of $N$ particles $\underline{x}(t)=\big(x_1(t),x_2(t),\ldots,x_N(t)\big)$ described by the mean-field model with Markovian switching
$$
    dx_i(t) = b\Big(x_i(t),\mu_N(t),\alpha(t_-)\Big)dt + \sigma\Big(x_i(t),\mu_N(t),\alpha(t_-)\Big) dw_i(t),\quad 1\le i\le N,
$$
where
$$
\mu_N(t)=\delta_{\underline{x}(t)}={1\over N}\sum_{j=1}^N\delta_{x_j(t)}\in\mathscr{M}_1.
$$
It is clear that this is a special case of the $N$-particle system given by \eqref{Eq-SDE-MSwit-a} with
\begin{equation}\label{Def-b_i,sigma_i}
\underline{b}_i\big(t,\underline{x},i_0\big)=b\big(x_i,\delta_{\underline{x}},i_0\big),\quad \underline{\sigma}_i\big(t,\underline{x},i_0\big)=\sigma\big(x_i,\delta_{\underline{x}},i_0\big),
\end{equation}
for $\big(t,\underline{x},i_0\big)\in[0,T]\times(\mathbb{R}^d)^N\times\mathbb{S}$.

 Note that $|\underline{x}|^2=\sum_{i=1}^N|x_i|^2$ implies $\big\langle \delta_{\underline{x}},\varphi\big\rangle={1\over N}\sum_{i=1}^N|x_i|\le{1\over\sqrt N}|\underline{x}|$ for any $\underline{x}\in(\mathbb{R}^d)^N$  and that
 $$
 \big\|\delta_{\underline{x}}-\delta_{\underline{y}}\big\|_{BL}\le {C\over N}\big|\underline{x}-\underline{y}\big|,\quad \forall\,\,\underline{x},\underline{y}\in(\mathbb{R}^d)^N.
 $$
Under Assumption (A), for $b(\cdot,\cdot,\cdot)$ and $\sigma(\cdot,\cdot,\cdot)$,
one can easily prove that the functions $\underline{b}_i$ and $\underline{\sigma}_i$, $1\le i\le N$, defined above satisfy the Lipschitz and linear growth conditions \eqref{Cond-SDE-Lips} and \eqref{Cond-SDE-Linear}. This implies that system \eqref{def-dynamics} has a unique solution.
The following lemma reveals the moment boundedness of the system. In order to keep the continuity of the presentation, its proof is relegated
to the Appendix.

\begin{lem}\label{lem-moments} Assume {\rm(A1)}, {\rm(A2)}, and that $\sup_{N\ge1}\EE\big\langle\mu_N(0),\psi\big\rangle<\infty$ where $\psi(x)=|x|^2$ for $x\in\mathbb{R}^d$. Then for positive numbers $T$ and $p$, $p\le1$, there is a constant $C$ independent of $N$ such that
\begin{equation}\label{moment-inq}
\sup_{0\le t\le T}\EE\bigg(\big\langle \mu_N(t),\psi\big\rangle+1\bigg)^p\le C,
\end{equation}
and for $0\le s\le t\le T$,
\begin{equation}\label{cond-moment-inq}
e^{-C(t-s)}\Big(\big\langle \mu_N(s),\psi\big\rangle+1\Big)^p\le \EE\Big[\Big(\big\langle \mu_N(t),\psi\big\rangle+1\Big)^p\Big|{\mathcal F}^{N,\alpha}_s\Big]\le e^{C(t-s)}\Big(\big\langle \mu_N(s),\psi\big\rangle+1\Big)^p.
\end{equation}
\end{lem}

For $f(\cdot,i_0)\in C^2_c\big(\mathbb{R}^d\big)$, $i_0\in\mathbb{S}$, and $\big(x,\mu,i_0\big)\in\mathbb{R}^d\times\mathscr{M}_1\times\mathbb{S}$ denote the operator
\begin{align}
\mathcal{L}(\mu)f\big(x,i_0\big)&=b'\big(x,\mu,i_0\big)\nabla_{x}f\big(x,i_0\big)+{1\over2}\Big(a\big(x,\mu,i_0\big)\nabla_{x}\Big)'\nabla_{x}f\big(x,i_0\big)\notag\\
&\quad+\sum_{j_0\in\mathbb{S}}q_{i_0j_0}\Big(f(x,j_0)-f(x,i_0)\Big),\label{Def-Op-MFSDE-MSw}
\end{align}
where
$$
a\big(x,\mu,i_0\big)=\sigma\big(x,\mu,i_0\big)\sigma'\big(x,\mu,i_0\big)\in\mathbb{R}^{d\times d}.
$$

Let $F(\cdot,\cdot,\cdot):\mathbb{R}\times(\mathbb{R}^d)^N\times\mathbb{S}\to\mathbb{R}$ be a function such that
$$
F\big(t,\underline{x},i_0\big)=\big\langle\delta_{\underline{x}},f(t,i_0)\big\rangle={1\over N}\sum_{i=1}^Nf\big(t,x_i,i_0\big),
$$
for some functions $f(\cdot,\cdot,i_0)\in C^{1,2}\big([0,T]\times\mathbb{R}^d\big)$, $i_0\in\mathbb{S}$. For $\underline{b}_i$ and $\underline{\sigma}_i$ defined as in \eqref{Def-b_i,sigma_i},  and $\mathcal{L}_N$ and $\mathcal{L}$ defined as in \eqref{Def-Op-SDE-MSw} and \eqref{Def-Op-MFSDE-MSw}, respectively,  we have
\begin{equation}
\bigg({\partial\over\partial t}+\mathcal{L}_N\bigg)F\big(t,\underline{x},i_0\big)={1\over N}\sum_{i=1}^N\bigg({\partial\over\partial t}+\mathcal{L}\big(\delta_{\underline{x}}\big)\bigg)f\big(t,x_i,i_0\big)=\bigg\langle \delta_{\underline{x}}, \bigg({\partial\over\partial t}+\mathcal{L}\big(\delta_{\underline{x}}\big)\bigg)f\big(t,i_0\big)\bigg\rangle.
\end{equation}

For each element $\varsigma \in D_f\big([0,T],\mathbb{S}\big)$ that represents a sample path of the switching process $\alpha(t)$, $0\le t\le T$, we denote the corresponding sample path of the associated martingale by a similar way to \eqref{Def-Quad-Var} and \eqref{Def-Mart of MC} as follow
\begin{equation}\label{Def-sample Mart of MC}
M_{i_0j_0}^{\varsigma}(t)=\big[M_{i_0j_0}^{\varsigma}\big](t)- \big\langle M_{i_0j_0}^{\varsigma}\big\rangle(t),
\end{equation}
where
$$\big[M_{i_0j_0}^{\varsigma}\big](t)=\sum_{0\le s\le t}{1\!\!1}\big(\varsigma(s_-)=i_0\big){1\!\!1}\big(\varsigma(s)=j_0\big),\quad \big\langle M_{i_0j_0}^{\varsigma}\big\rangle(t) =\int_0^tq_{i_0j_0}{1\!\!1}\big(\varsigma(s_-)=i_0\big)ds.
$$
Since the sample paths of $\alpha(\cdot)$ are in $D_f\big([0,T],\mathbb{S}\big)$, for simplicity, we define $M_{i_0j_0}^{\varsigma}(t)=\big[M_{i_0j_0}^{\varsigma}\big](t)=\big\langle M_{i_0j_0}^{\varsigma}\big\rangle(t)=0$ for $\varsigma\in D\big([0,T],\mathbb{S}\big)\backslash D_f\big([0,T],\mathbb{S}\big)$ and $0\le t\le T$.

For $\big(\underline{x}(\cdot),\alpha(\cdot)\big) \in D\big([0,T],(\mathbb{R}^d)^N\times\mathbb{S}\big)$,
 we define the mapping $e_N$ by
$$
e_N:
\big(\underline{x}(\cdot),\alpha(\cdot)\big)\mapsto\big(\delta_{\underline{x}(\cdot)},\alpha(\cdot)\big)\in D\big([0,T],\mathscr{M}_1\times\mathbb{S}\big).
$$
Denote by $\PP_N$ the induced probability measure of the system $\big(x_1(\cdot),x_2(\cdot),\ldots,x_N(\cdot),\alpha(\cdot)\big)$, the solution to \eqref{def-dynamics}, on $D\big([0,T],(\mathbb{R}^d)^N\times\mathbb{S}\big)$. It follows that
$$
\mathscr{P}_N=\PP_N\circ e_N^{-1},\quad \mu_N(\cdot)=\delta_{\underline{x}(\cdot)}.
$$
We have the following lemma.

\begin{lem}\label{lem-martingale} Under the assumption of \lemref{lem-moments}, the following statements hold:
\begin{itemize}
\item[{\rm(i)}] For $f(\cdot,\cdot,i_0)\in C^{1,2}\big([0,T]\times\mathbb{R}^d\big)$, $i_0\in\mathbb{S}$, and $\big(\eta,\varsigma\big)\in D\big([0,T],\mathscr{M}_1\times\mathbb{S}\big)$ denote
\begin{align}
M_f(t)=\big\langle\eta(t),f\big(t,\cdot,\varsigma(t)\big)\big\rangle
&-\big\langle\eta(0),f\big(0,\cdot,\varsigma(0)\big)\big\rangle\notag\\
& -\int_0^t\bigg\langle\eta(s),\bigg({\partial \over\partial s}+\mathcal{L}\big(\eta(s)\big)\bigg)f\big(s,\cdot,\varsigma(s_-)\big)\bigg\rangle ds\notag\\
&-\sum_{i_0,j_0\in\mathbb{S}}\int_0^t\Big\langle\eta(s),f\big(s,\cdot,j_0\big)-f\big(s,\cdot,i_0\big)\Big\rangle dM^{\varsigma}_{i_0j_0}(s).\label{Def-M^f_t}
\end{align}
Then $M_f(t)$ is a continuous $\mathscr{P}_N-$martingale.
\item[{\rm(ii)}] For $f(\cdot,\cdot,i_0), g(\cdot,\cdot,i_0)\in C^{1,2}\big([0,T]\times\mathbb{R}^d\big)$, $i_0\in\mathbb{S}$, the quadratic variational process of the $\mathscr{P}_N-$martingales $M_f$ and $M_g$ has the form
\begin{align*}
[M_f,M_g](t)&=[M_f,M_g]_N(t)\\
&={1\over N}\int_0^t\bigg\langle\eta(s),\Big(a\big(\cdot,\eta(s),\varsigma(s_-)\big)\nabla f(s,\cdot,\varsigma(s_-)),\nabla g(s,\cdot,\varsigma(s_-))\Big)\bigg\rangle ds,
\end{align*}
for $0\le t\le T$, where $\nabla$ denotes the gradient with respect to the space variables and $(\cdot,\cdot)$ is the inner product in $\mathbb{R}^d$.
\end{itemize}
\end{lem}

\para{Proof.}
Let $F(\cdot,\cdot,\cdot), G(\cdot,\cdot,\cdot):\mathbb{R}\times(\mathbb{R}^d)^N\times\mathbb{S}\to\mathbb{R}$ satisfy
$$
F\big(t,\underline{x},i_0\big)=\big\langle\delta_{\underline{x}},f\big(t,\cdot,i_0\big)\big\rangle,\quad G\big(t,\underline{x},i_0\big)=\big\langle\delta_{\underline{x}},g\big(t,\cdot,i_0\big)\big\rangle,
$$
for functions $f(\cdot,\cdot,i_0), g(\cdot,\cdot,i_0)\in C^{1,2}\big([0,T]\times\mathbb{R}^d\big)$, $i_0\in\mathbb{S}$. We put
\begin{align}
M^F(t)= &\,\, F\big(t,\underline{x}(t),\alpha(t)\big)-F\big(0,\underline{x}(0),\alpha(0)\big)-\int_0^t\bigg({\partial\over\partial s}+\mathcal{L}_N\bigg)F\big(s,\underline{x}(s),\alpha(s_-)\big)ds\notag\\
&-\sum_{i_0,j_0\in\mathbb{S}}\int_0^t\Big(F\big(t,\underline{x}(t),j_0\big)-F\big(t,\underline{x}(t),i_0\big)\Big)dM_{i_0j_0}(s),\quad 0\le t\le T,\notag
\end{align}
and define $M^G(t)$ with $F$ replaced by $G$ and $0\le t\le T$.  In view of
the It\^o formula \eqref{Eq-Ito-SDEMSw},
\begin{align}
M^F(t)&=\sum_{i=1}^N\int_0^t\Big\langle\nabla_{x_i}F\big(s,\underline{x}(s),\alpha(s_-)\big),\underline{\sigma}_i\big(s,\underline{x}(s),\alpha(s_-)\big)dw_i(s)\Big\rangle\notag\\
&={1\over N}\sum_{i=1}^N\int_0^t\Big\langle\nabla_{x}f\big(s,x_i(s),\alpha(s_-)\big),\sigma\big(x_i(s),\delta_{\underline{x}(s)},\alpha(s_-)\big)dw_i(s),\Big\rangle\label{Def-M^F}
\end{align}
is a continuous $\PP_N$-martingle. Since the Brownian motions $w_1(\cdot),w_2(\cdot),\ldots,w_N(\cdot)$ are independent, from \eqref{Def-M^F}, we obtain
\begin{align}
&[M^F,M^G](t)\notag\\
&={1\over N^2}\sum_{i=1}^N\int_0^t\Big(a\big(x_i(s),\delta_{\underline{x}(s)},\alpha(s_-)\big)\nabla f\big(s,x_i(s),\alpha(s_-)\big),\nabla g\big(s,x_i(s),\alpha(s_-)\big)\Big)ds\notag\\
&={1\over N}\int_0^t\bigg\langle\delta_{\underline{x}(s)},\Big(a\big(\cdot,\delta_{\underline{x}(s)},\alpha(s_-)\big)\nabla f\big(s,\cdot,\alpha(s_-)\big),\nabla g\big(s,\cdot,\alpha(s_-)\big)\Big)\bigg\rangle ds.\label{Eq-[M^F,M^G]}
\end{align}
One can easily verify the following identities
$$
M^F=M_f\circ e_N,\quad M^G=M_g\circ e_N, \quad [M^F,M^G]=[M_f,M_g]\circ e_N.
$$
Since $\mathscr{P}_N=\PP_N\circ e_N^{-1}$, a combination of the above facts implies the assertions (i) and (ii).\qed

\section{Law of Large Numbers for Mean-Field Models with Markovian Switching}\label{sec:lln}

In this section, we present the proof of \thmref{thm-main1}, one of the main result of the paper, establishing the law of large numbers for the mean-field systems with Markovian switching. ​We use the martingale approach. The weak compactness of the sequence $\{(\mu_N, \alpha_N)\}_{N\ge1}$ is established in Section 4.1. Its limit is characterized in Section 4.2.

\subsection{Weak Compactness of $\big\{\big(\mu_N(\cdot),\alpha(\cdot)\big)\big\}_{N\ge1}$}

\begin{lem}\label{lem-compact-containment} Under the assumptions of \thmref{thm-main1}, for each $\delta>0$ there exists a compact set $K_\delta$ in $\big(\mathscr{M}_1,\|\cdot\|_{BL}\big)$ such that
$$\inf_{N\ge1}\PP\Big(\mu_N(t)\in K^\delta_\delta\,\text{ for all }\,  0\le t\le T\Big)\ge1-\delta,$$
where $K^\delta_\delta=\big\{\mu\in\mathscr{M}_1:\inf_{\eta\in K_\delta}\|\mu-\eta\|_{BL}<\delta\big\}$.
\end{lem}

\para{Proof.} For each $\lambda>0$, denote $B_\lambda^c=\big\{x\in\mathbb{R}^d: |x|>\lambda\big\}$ and $H_\lambda=\big\{\mu\in\mathscr{M}_1:\mu(B^c_\lambda)=0\big\}$. Because $\big(\mathscr{M}_1,\|\cdot\|_{BL}\big)$ is topologically equivalent to $\mathcal{P}(\mathbb{R}^d)$ (see Theorem 12 \cite{Dudley66}), by Prohorov theorem, $H_\lambda$ is relatively compact in $\big(\mathscr{M}_1,\|\cdot\|_{BL}\big)$. For $0\le t\le T$,
\begin{align}
\Big\{\mu_N(t)\notin H^\delta_\lambda\Big\}
&\subset\Big\{\forall\, \eta\in H_\lambda:\big\|\eta-\mu_N(t)\big\|_{TV}\ge\delta\Big\}\notag\\
&=\bigg\{\forall\, \eta\in H_\lambda:\Big\|\Big(\eta-\mu_N(t)\Big)_{|_{B_\lambda}}\Big\|_{TV}+\Big\|\mu_N(t)_{|_{B^c_\lambda}}\Big\|_{TV}\ge\delta\bigg\}\notag\\
&\subset\bigg\{\Big\|\mu_N(t)_{|_{B^c_\lambda}}\Big\|_{TV}\ge{\delta\over2}\bigg\}\notag\\
&\subset\bigg\{\big\langle\mu_N(t),\psi\big\rangle\ge{\lambda^2\delta\over2}\bigg\}.\label{inq-compact}
\end{align}
Let $C$ be the constant given in \eqref{cond-moment-inq}. It follows from \eqref{inq-compact} and \eqref{cond-moment-inq} that
\begin{align*}
\PP\Big(\mu_N(t)\notin H^\delta_\lambda\,\,\text{ for some }\,\,0\le t\le T\Big)
&\le \PP\bigg(\sup_{0\le t\le T}e^{Ct}\big\langle\mu_N(t),\psi\big\rangle\ge{\lambda^2\delta\over2}\bigg)\\
&\le {2\over \lambda^2\delta}\EE\bigg[e^{CT}\big\langle\mu_N(T),\psi\big\rangle\bigg]\\
&\le {2\over \lambda^2\delta}\EE\bigg[e^{2CT}\big\langle\mu_N(0),\psi\big\rangle\bigg]\\
&={C\over\lambda^2\delta}.
\end{align*}
For a fixed $\delta>0$, we can choose $\lambda=\lambda(\delta)$ large enough such that ${C\over\lambda^2\delta}\le\delta$. Take $K_\delta=\bar H_{\lambda(\delta)}$ which is compact in $\big(\mathscr{M}_1,\|\cdot\|_{BL}\big)$. For all $N\ge1$ we have
\begin{align*}
\PP\Big(\mu_N(t)\in K^\delta_\delta\text{ for all }\,  0\le t\le T\Big)\ge1-\PP\Big(\mu_N(t)\notin H^\delta_{\lambda(\delta)}\,\,\text{ for some }\,\,0\le t\le T\Big)\ge 1-\delta.
\end{align*}
This completes the proof. \qed

\begin{lem}\label{lem-weak-compactness} Under the assumptions of \thmref{thm-main1}, for each positive integer $N$ and $\delta>0$, there exists a random variable  $\gamma_N(\delta)\ge0$ such that
$$
\EE\Big(\big\|\mu_N(t+\delta)-\mu_N(t)\big\|_{BL}\Big|\mathcal{F}^{N,\alpha}_t\Big)\le \EE\Big(\gamma_N(\delta)\Big|\mathcal{F}^{N,\alpha}_t\Big)\quad\text{a.s.}
$$
 for all $0\le t\le T-\delta$. Furthermore,
$$
\lim_{\delta\to0}\limsup_{N\to\infty}\EE\gamma_N(\delta)=0.
$$
\end{lem}

\para{Proof.} By the definition of the norm $\|\cdot\|_{BL}$ and Cauchy-Schwarz inequality,
\begin{align*}
\big\|\mu_N(t+\delta)-\mu_N(t)\big\|_{BL}^2
&\le{1\over N}\sum_{i=1}^N\big|x_i(t+\delta)-x_i(t)\big|^2,
\end{align*}
for any integer $N$, and real numbers $t,\delta$ satisfying $0\le t,t+\delta\le T$. Therefore, by the Dynkin formula,
\begin{align}
&\EE\Big(\big\|\mu_N(t+\delta)-\mu_N(t)\big\|^2_{BL}\Big|\mathcal{F}^{N,\alpha}_t\Big)\notag\\
&\le \EE\bigg({1\over N}\sum_{i=1}^N\big|x_i(t+\delta)-x_i(t)\big|^2\bigg|\mathcal{F}^{N,\alpha}_t\bigg)\notag\\
&= {1\over N}\sum_{i=1}^N\EE\bigg\{\int_t^{t+\delta}\bigg[2b'\big(x_i(s),\mu_N(s),\alpha(s_-)\big)\Big(x_i(s)-x_i(t)\Big)\notag\\
&\hspace{6cm}+\text{trace}\Big(a\big(x_i(s),\mu_N(s),\alpha(s_-)\big)\Big)\bigg]ds\bigg|\mathcal{F}^{N,\alpha}_t\bigg\}\notag\\
&\le {1\over N}\sum_{i=1}^N\EE\bigg\{\int_t^{t+\delta}\bigg[2\Big(\big|x_i(s)\big|+\big\langle\mu_N(s),\varphi\big\rangle+1\Big)\Big(\big|x_i(s)\big|+\big|x_i(t)\big|\Big)\notag\\
&\hspace{6cm}
+3\Big(\big|x_i(s)\big|^2+\big\langle\mu_N(s),\psi\big\rangle+1\Big)\bigg]ds\bigg|\mathcal{F}^{N,\alpha}_t\bigg\}\notag\\
&\le C\EE\bigg[\int_t^{t+\delta}\Big(\big\langle\mu_N(s),\psi\big\rangle+\big\langle\mu_N(t),\psi\big\rangle+1\Big)ds\bigg|\mathcal{F}^{N,\alpha}_t\bigg].\label{est-weak-comp}
\end{align}
It follows from the right-hand side and then the left-hand side of \eqref{cond-moment-inq} that for $s\ge t$,
$$
\EE\Big(\langle\mu_N(s),\psi\rangle+1\Big|\mathcal{F}^{N,\alpha}_t\Big)\le C\Big(\langle\mu_N(t),\psi\rangle+1\Big)\le C\EE\Big(\langle\mu_N(T),\psi\rangle+1\Big|\mathcal{F}^{N,\alpha}_t\Big).
$$
Thus, \eqref{est-weak-comp} implies that
\begin{align}
\EE\Big(\big\|\mu_N(t+\delta)-\mu_N(t)\big\|^2_{BL}\Big|\mathcal{F}^{N,\alpha}_t\Big)&\le C\int_t^{t+\delta}\Big(\big\langle\mu_N(t),\psi\big\rangle+1\Big)ds\notag\\
&\le  C\int_t^{t+\delta}\EE\Big(\big\langle\mu_N(T),\psi\big\rangle+1\Big|\mathcal{F}^{N,\alpha}_t\Big)ds\notag\\
&=C\delta \EE\Big(\big\langle\mu_N(T),\psi\big\rangle+1\Big|\mathcal{F}^{N,\alpha}_t\Big).\notag
\end{align}
As a consequence, by Cauchy-Schwarz inequality,
\begin{align}
\EE\Big(\big\|\mu_N(t+\delta)-\mu_N(t)\big\|_{BL}\Big|\mathcal{F}^{N,\alpha}_t\Big)&\le \Big[\EE\Big(\big\|\mu_N(t+\delta)-\mu_N(t)\big\|^2_{BL}\Big|\mathcal{F}^{N,\alpha}_t\Big)\Big]^{1/2}\notag\\
&\le {1\over2}\Big[\sqrt{\delta}+C\sqrt{\delta} \EE\Big(\big\langle\mu_N(T),\psi\big\rangle+1\Big|\mathcal{F}^{N,\alpha}_t\Big)\Big].\label{inq-mu-diff}
\end{align}
 This inequality and \lemref{lem-moments} conclude the proof by taking
 $$
 \gamma_N(\delta)=\sqrt{\delta}C\Big(\big\langle\mu_N(T),
 \psi\big\rangle+1\Big).
 $$
 The proof is complete.
\qed

According to Lemmas \ref{lem-compact-containment} and \ref{lem-weak-compactness} we obtain the following Proposition.

\begin{prop}
The sequence $\big\{\big(\mu_N(\cdot),\alpha(\cdot)\big), N\ge1\big\}$ is weakly compact in the topology of weak convergence of probability measure on $D\big([0,T],\mathscr M_1\times\mathbb S\big)$.
\end{prop}

\para{Proof.}
Let $K_\delta$ be a compact subset of $\big(\mathscr{M}_1,\|\cdot\|_{BL}\big)$ as in \lemref{lem-compact-containment} and denote $L_\delta=K_\delta\times\mathbb{S}$. By the compactness of the space $(\mathbb{S},d_\mathbb{S})$, $L_\delta$ is also a compact set in $\big(\mathscr{M}_1\times\mathbb{S},d\big)$ and $L_\delta^\delta=K_\delta^\delta\times\mathbb{S}$. In view of \lemref{lem-compact-containment}, we have
\begin{equation}\label{Eq-compact-containment}
\inf_{N\ge1}\PP\Big(\big(\mu_N(t),\alpha(t)\big)\in L_\delta^\delta\,\,\,\,\text{for all}\,\,\,\,0\le t\le T\Big)\ge1-\delta.
\end{equation}
Since $\alpha(\cdot)$ is a Markov chain, there is a constant $C$ such that $\EE\big(d_\mathbb{S}(\alpha(t+\delta),\alpha(t))\big|\mathcal{F}^\alpha_t\big)\le C\delta$ for $\delta>0$ and $0\le t\le T$. Therefore, it follows from the definition of metric $d$ and \lemref{lem-weak-compactness} that for each integer $N$ and positive number $\delta$ there exists a random variable $\gamma_N(\delta)$ such that
\begin{align}
&\EE\bigg[d\Big(\big(\mu_N(t),\alpha(t)\big),\big(\mu_N(t+\delta),\alpha(t+\delta)\big)\Big)\bigg|\mathcal{F}^{N,\alpha}_t\bigg]\notag\\
&=\EE\bigg[\big\|\mu_N(t+\delta)-\mu_N(t)\big\|_{BL}+d_{\mathbb{S}}\big(\alpha(t+\delta),\alpha(t)\big)\bigg|\mathcal{F}^{N,\alpha}_t\bigg]\notag\\
&\le \EE\Big(\gamma_N(\delta)+C\delta\Big|\mathcal{F}^{N,\alpha}_t\Big), \label{Eq-weak-compactness}
\end{align}
where $\gamma_N(\delta)$ defined in \lemref{lem-weak-compactness} satisfies
\begin{equation}\label{Eq-weak-compactness1}
 \lim_{\delta\to0} \limsup_{N\to\infty}\EE
 \Big(\gamma_N(\delta)+C\delta \Big)=0.
\end{equation}
Combining \eqref{Eq-compact-containment}, \eqref{Eq-weak-compactness}, and \eqref{Eq-weak-compactness1}, the Proposition follows by virtue of   \cite[Theorem 3.8.6]{EthierKurtz86}. \qed

\subsection{Characterization of Limit}
Next, we proceed to characterize
$\big(\mu_\alpha(\cdot),\alpha(\cdot)\big)$, the limit of  the sequence $\big\{(\mu_N(\cdot),\alpha(\cdot))\big\}_{ N\ge1}$.
 We have the following lemma.

\begin{thm}\label{thm-limit-characterization} Assume {\rm(A1)}, {\rm(A2)}, and that $\sup_{N\ge1}\EE\big\langle\mu_N(0),\psi\big\rangle<\infty$. Denote by $\mathscr{P}$ the limit of an arbitrary weakly convergent subsequence of $\mathscr{P}_N$. Then for $\mathscr{P}$-almost all $\big(\eta(\cdot),\varsigma(\cdot)\big)\in D\big([0,T],\mathscr{M}_1\times\mathbb{S}\big)$,
\begin{align}
\big\langle\eta(t),f\big(\cdot,\varsigma(t)\big)\big\rangle=\big\langle\eta(0),f\big(\cdot,\varsigma(0)\big)\big\rangle&+\int_0^t\Big\langle\eta(s),\mathcal{L}\big(\eta(s)\big)f\big(\cdot,\varsigma(s_-)\big)\Big\rangle ds\notag\\
&+\sum_{i_0,j_0\in\mathbb{S}}\int_0^t\Big\langle\eta(s),f\big(\cdot,j_0\big)-f\big(\cdot,i_0\big)\Big\rangle dM_{i_0j_0}^\varsigma(s), \label{eq-stochastic-MV}
\end{align}
holds for all test functions $f(\cdot,i_0)\in C^{2}_c\big(\mathbb{R}^d\big)$, $i_0\in\mathbb{S}$, and $0\le t\le T$.
\end{thm}

\para{Proof.}
 Let $\mathscr{P}_{N_k}$ be a subsequence of $\mathscr{P}_{N}$ that  converges weakly  to a probability measure $\mathscr{P}$ on $D\big([0,T],\mathscr{M}_1\times\mathbb{S}\big)$ as $k\to\infty$. It suffices to prove that \eqref{eq-stochastic-MV} holds $\mathscr{P}-$almost surely for each test function $f(\cdot,i_0)\in C^{2}_c\big(\mathbb{R}^d\big)$, $i_0\in\mathbb{S}$. Note from Lemma \ref{lem-martingale} that for each $N$,
\begin{align}
 M_f(t)=\big\langle\eta(t),f\big(\cdot,\varsigma(t)\big)\big\rangle&-\big\langle\eta(0),f\big(\cdot,\varsigma(0)\big)\big\rangle-\int_0^t\Big\langle\eta(s),\mathcal{L}(\eta)f\big(\cdot,\varsigma(s)\big)\Big\rangle ds\notag\\
&-\sum_{i_0,j_0\in\mathbb{S}}\int_0^t\Big\langle\eta(s),f\big(\cdot,j_0\big)-f\big(\cdot,i_0\big)\Big\rangle dM^\varsigma_{i_0j_0}(s),\quad 0\le t\le T, \label{eq-mart}
\end{align}
is a continuous square integrable $\mathscr{P}_{N}$-martingale with quadratic variational process
$$
[M_f]_N(t)={1\over N}\int_0^t\bigg\langle\eta(s),\Big(a\big(\cdot,\eta(s),\varsigma(s_-)\big)\nabla f(\cdot,\varsigma(s_-)),\nabla f(\cdot,\varsigma(s_-))\Big)\bigg\rangle ds,\quad 0\le t\le T.
$$
By virtue of Assumption (A2) and Lemma \ref{lem-moments}, there exists a constant $C$ independent of $N$ such that
\begin{equation}\label{est-quad-var}
\sup_{0\le t\le T}\EE_{\mathscr{P}_{N}}[M_f]_{N}(t)\le {CT\over N}.
\end{equation}
Note that for any $N$, and $0\le t\le T$, $(\mu_N(\cdot),\alpha(\cdot))$ concentrates on the set $\mathscr{S}$ denoted by
\begin{align*}
\mathscr{S}= &  \bigg\{\big(\eta,\varsigma\big)\in C\big([0,T],\mathscr{M}_1\big)\times D_f\big([0,T],\mathbb{S}\big): \varsigma(t)=\varsigma(t_-) \bigg\},
\end{align*}
 with probability $1$, i.e., $\mathscr{P}_N\big(\mathscr{S}\big)=1$ for any $N\ge1$. For $n\ge1$, denote
 $$
 \mathscr{S}_n=\Big\{\big(\eta,\varsigma\big)\in\mathscr{S}: \varsigma \text{ has no jump in $\big(t-t/(2n),t\big]$}\Big\}.
 $$
 It is clear that $\mathscr{S}_n\subset \mathscr{S}_{n+1}$ for every $n\ge1$ and $\mathscr{S}=\cup_{n\ge1}\mathscr{S}_n$.  Since $\mathbb{S}$ is a discrete set and $f\in C^2_c(\mathbb{R}^d\times\mathbb{S})$, the set $\big\{(\eta,\varsigma)\in \mathscr{S}_n :|M_f(t)|\le\delta\big\}$ is closed in $D([0,T],\mathscr{M}_1\times\mathbb{S})$ for each $n\ge1$. Thus, by Portmanteau Theorem, we have
\begin{align}
\mathscr{P}\Big(\big|M_f(t)\big|\le\delta\Big)
&\ge\mathscr{P}\bigg((\eta,\varsigma)\in\mathscr{S}:|M_f(t)|\le\delta \bigg)\notag\\
&=\lim_{n\to\infty}\mathscr{P}\bigg((\eta,\varsigma)\in\mathscr{S}_n:|M_f(t)|\le\delta\bigg)\notag\\
&\ge \lim_{n\to\infty}\limsup_{k\to\infty}\mathscr{P}_{N_k}\bigg((\eta,\varsigma)\in\mathscr{S}_n:|M_f(t)|\le\delta\bigg)\notag\\
&\ge\lim_{n\to\infty}\bigg[\mathscr{P}_{N_k}\Big(\mathscr{S}_n\Big)-\limsup_{k\to\infty}\mathscr{P}_{N_k}\Big(\big|M_f(t)\big|>\delta\Big)\bigg]\notag\\
&=1-\limsup_{k\to\infty}\mathscr{P}_{N_k}\Big(\big|M_f(t)\big|>\delta\Big).\label{est-ball-quad-var}
\end{align}
By Doob's submartingale inequality and \eqref{est-quad-var},
\begin{equation}\label{est-tail-quad-var}
\mathscr{P}_{N_k}\big(|M_f(t)|>\delta\big)\le \delta^{-2}\EE_{\mathscr{P}_{N_k}}[M_f]_{N_k}(t)\le {CT\over\delta^2N_k}\to0 \quad \text{as}\quad k\to\infty.
\end{equation}
Combining \eqref{est-ball-quad-var} and \eqref{est-tail-quad-var} yields $\mathscr{P}\big(|M_f(t)|=0\big)=1$, which implies $M_f(t)=0$ $\mathscr{P}$-a.s. The theorem therefore follows from \eqref{eq-mart}.\qed

To proceed, we need  a result from \cite{Oelschlager84}. Assume the following conditions hold.
\begin{itemize}
\item[(B1)] $\hat b(\cdot,\cdot):\mathbb{R}^d\times\mathscr M_1\to\mathbb{R}^d$ and $\hat \sigma(\cdot,\cdot):\mathbb{R}^d\times\mathscr M_1\to\mathbb{R}^{d\times d}$ are Lipschitz continuous in that, there is a constant $L$ such that
    $$
    \big|\hat b\big(x,\mu\big)-\hat b\big(y,\eta\big)\big|+\big|\hat \sigma\big(x,\mu\big)-\hat \sigma\big(y,\eta\big)\big|\le L\big(\big|x-y\big|+\big\|\mu-\eta\big\|_{BL}\big),
    $$
    for all $x,y\in\mathbb{R}^d$ and $\mu,\eta\in \mathscr M_1$.

\item[(B2)]  The $\mathbb{R}^d$-valued function $\hat b(\cdot,\cdot)$ satisfies
 $$
 \big|\hat b\big(x,\mu\big)\big|\le C\big(1+\big|x\big|+\big\langle\mu,\varphi\big\rangle\big),\quad  (x,\mu)\in\mathbb{R}^d\times\mathscr{M}_1,
 $$
for some constant $C$ and the $\mathbb{R}^{d\times d}$-valued function $\hat \sigma(\cdot,\cdot)$ is bounded.
\end{itemize}

For $\mu\in \mathscr M_1$ and $f\in C^2_c(\mathbb{R}^d)$ denote
$$
\hat{\mathcal{L}}(\mu)f(x)=\hat b'(x,\mu)\nabla_xf(x)+{1\over2}\big(\hat a(x,\mu)\nabla_x\big)'\nabla_xf(x),\quad x\in\mathbb{R}^d.
$$
As a consequence of Lemma 9 and equation (8.2) in \cite{Oelschlager84}, we have the following theorem.

\begin{thm}\label{thm-Oelschlager}
Assume {\rm(B1)} and {\rm(B2)}.  Then the equation
\begin{equation}\label{eq-SMV-O}
\big\langle \mu(t),f\big\rangle =\big\langle \mu_0,f\big\rangle + \int_0^t\Big\langle \mu(s),\hat{\mathcal{L}}\big(\mu(s)\big)f\Big\rangle ds,\quad 0\le t \le T, f\in C^2_c(\mathbb{R}^d),
\end{equation}
has a unique solution $\mu(t)=\mathscr{L}(z(t))$ in $D\big([0,T],\mathscr{M}_1\big)$ that is the distribution of the unique solution of
$$
dz(t)=\hat b\Big(z(t),\mathscr L\big(z(t)\big)\Big)dt+\hat \sigma\Big(z(t),\mathscr L\big(z(t)\big)\Big)d\tilde w(t),\quad \mathscr L(z_0)=\mu_0,
$$
and $\tilde w(\cdot)$ is a standard Brownian motion.
\end{thm}

We are now in a position to present a result on stochastic McKean-Vlasov equations with Markovian switching.

\begin{thm}\label{thm-SMV}
Assume {\rm(A1)} and {\rm(A2)}.  Then the system of integral equations
\begin{align}
\big\langle\mu(t),f(\cdot,\alpha(t))\big\rangle=  \big\langle\mu_0,f(\cdot,\alpha(0))\big\rangle&+\int_0^t\Big\langle\mu(s),\mathcal{L}\big(\mu(s)\big)f\big(\cdot,\alpha(s_-)\big)\Big\rangle ds\notag\\
&+\sum_{i_0,j_0\in\mathbb{S}}\int_0^t\Big\langle\mu(s),f(\cdot,j_0)-f(\cdot,i_0)\Big\rangle dM_{i_0j_0}(s), \label{eq-perturb-SMV}
\end{align}
where $0\le t\le T$ and $f(\cdot,i_0)\in C^2_c(\mathbb{R}^d)$ for each $i_0\in \mathbb{S}$, has a unique solution in $D\big([0,T],\mathscr{M}_1\big)$. Moreover, this solution equals $\mathscr{L}\big(y(t)\big|\mathcal{F}^\alpha_{t_-}\big)$ for all $0\le t\le T$,
 where $y(t)$ is the unique solution of
$$
\begin{cases}
dy(t)&=b\Big(y(t),\mu_\alpha(t),\alpha(t_-)\Big)dt+\sigma\Big(y(t),\mu_\alpha(t),\alpha(t_-)\Big)d\tilde w(t),\quad \mathscr L(y(0))=\mu_0,\\
\mu_\alpha(t)&=\mathscr{L}\big(y(t)\big|\mathcal{F}^\alpha_{t_-}\big),
\end{cases}
$$
where $\tilde w(\cdot)$ is a standard Brownian motion independent of $\alpha(\cdot)$.
\end{thm}

Since the proof of this theorem is rather long, we give  a brief explanation of the main idea. First observe that \eqref{eq-SMV-O} is a special case of \eqref{eq-perturb-SMV} for the mean-field models with the usual diffusion (i.e., without switching process).  To proceed with the case involving Markovian switching, we use \thmref{thm-Oelschlager} to deal with  equation \eqref{eq-perturb-SMV} in the time intervals between the jumps of the Markov chain $\alpha(\cdot)$. We then consider \eqref{eq-perturb-SMV} at any jump time point  by ``gluing'' the solutions between jump times of the Markov chain and show that the solution obtained in this way indeed satisfies all the requirements.

\para{Proof of \thmref{thm-SMV}.} The proof is divided into several steps.

\underline{Step 1:}  Show that for each $0<r\le T$ and $\varsigma(\cdot)\in D_f\big([0,r],\mathbb{S}\big)$, there exists a unique solution $\eta(\cdot)\in D\big([0,r],\mathscr{M}_1\big)$ to the equation
\begin{align}
\big\langle\eta(t),f\big(\cdot,\varsigma(t)\big)\big\rangle=\big\langle\mu_0,f\big(\cdot,\varsigma(0)\big)\big\rangle&+\int_0^t\Big\langle\eta(s),\mathcal{L}(\eta)f\big(\cdot,\varsigma(s_-)\big)\Big\rangle ds\notag\\
&+\sum_{i_0,j_0\in\mathbb{S}}\int_0^t\Big\langle\eta(s),f(\cdot,j_0)-f(\cdot,i_0)\Big\rangle dM^\varsigma_{i_0j_0}(s), \label{eq-SMV-varsig}
\end{align}
where $0\le t\le r$ and $f(\cdot,i_0)\in C^2_c(\mathbb{R}^d)$ for each $i_0\in \mathbb{S}$.

Denote $t_0=0$, $t_{n+1}=\inf\{t>t_n:\varsigma(t)\ne\varsigma(t_-)\}$ and $\iota_n=\varsigma(t_n)\in\mathbb{S}$ for $n\ge0$.
For each $i_0\in\mathbb{S}$, $x\in\mathbb{R}^d$ and $\mu\in {\mathscr{M}}_1$, denote $\hat b_{i_0}(x,\mu)=b(x,\mu,i_0)$, $\hat \sigma_{i_0}(x,\mu)=\sigma(x,\mu,i_0)$, and $\hat a_{i_0}(x,\mu)=\hat a(x,\mu,i_0)$. In addition, for $\mu\in {\mathscr M}_1$, $f\in C^2_c(\mathbb{R}^d)$ denote
$$
\hat{\mathcal{L}}_{i_0}(\mu)f(x)=\hat b_{i_0}'(x,\mu)\nabla_xf(x)+{1\over2}\big(\hat a_{i_0}(x,\mu)\nabla_x\big)'\nabla_xf(x),\quad x\in\mathbb{R}^d.
$$
Then
\begin{equation}\label{eqn-L=L^}
\mathcal{L}(\mu)f(x,i_0)= \hat{\mathcal{L}}_{i_0}(\mu)f(x,i_0)+\sum_{j_0\in\mathbb{S}}q_{i_0j_0}\Big(f(x,j_0)-f(x,i_0)\Big).
\end{equation}

Next, we show that on each interval $[t_{k-1},t_k]$, $k\ge1$, a solution to \eqref{eq-SMV-varsig} satisfies  equation \eqref{eq-SMV-O} in \thmref{thm-Oelschlager} with the operator $\hat{\mathcal{L}}_{\iota_{k-1}}$. That is, for $t_{k-1}\le t\le t_k$,
\begin{equation}
\langle\eta(t),f(\cdot,\iota_{k-1})\rangle=\langle\eta(t_{k-1}),f(\cdot,\iota_{k-1})\rangle+\int_{t_{k-1}}^t\Big\langle\eta(s),\hat{\mathcal{L}}_{\iota_{k-1}}\big(\eta(s)\big)f(\cdot,\iota_{k-1})\Big\rangle ds.\label{eq-SMV-varsig-k}
\end{equation}
First, we consider the case $k=1$. For $t_0\le t<t_1$ we have $M^\varsigma_{i_0j_0}(t)=-\big\langle M^\varsigma_{\iota_0j_0}\big\rangle(t)=-q_{\iota_0j_0}t$ if $j_0\ne i_0=\iota_0$, and $M^\varsigma_{i_0j_0}(t)=0$ if $i_0\ne\iota_0$. Thus it follows from \eqref{eq-SMV-varsig} and \eqref{eqn-L=L^} that for any $f\in C^2_c(\mathbb{R}^d\times\mathbb{S})$ and $0<t<t_1$,
\begin{align}
&\big\langle\eta(t),f(\cdot,\iota_0)\big\rangle\notag\\
&= \big\langle\mu_0,f(\cdot,\iota_0)\big\rangle+\int_0^t\Big\langle\eta(s),\mathcal{L}\big(\eta(s)\big)f(\cdot,\iota_0)\Big\rangle ds   -\sum_{j_0\in\mathbb{S}}\int_0^t\Big\langle\eta(s),f(\cdot,j_0)-f(\cdot,\iota_0)\Big\rangle d\big\langle M^\varsigma_{\iota_0j_0}\big\rangle(s)\notag\\
&=\big\langle\mu_0,f(\cdot,\iota_0)\big\rangle+\int_0^t\Big\langle\eta(s),\hat{\mathcal{L}}_{\iota_0}\big(\eta(s)\big)f(\cdot,\iota_0)\Big\rangle ds.\label{eq-SMV-varsig-I1}
\end{align}
By taking $t=t_1$ in \eqref{eq-SMV-varsig} and noting that $\big\langle M^\varsigma_{\iota_0j_0}\big\rangle(t_1)=q_{\iota_0j_0}t_1$ for $j_0\ne\iota_0$ and $\big[ M^\varsigma_{\iota_0\iota_1}\big](t_1)=1$, we obtain
\begin{align}
&\big\langle\eta(t_1),f(\cdot,\iota_1)\big\rangle \notag\\
&= \big\langle\mu_0,f(\cdot,\iota_0)\big\rangle+\int_0^{t_1}\Big\langle\eta(s),\mathcal{L}\big(\eta(s)\big)f(\cdot,\iota_0)\Big\rangle ds   -\sum_{j_0\in\mathbb{S}}\int_0^{t_1}\Big\langle\eta(s),f(\cdot,j_0)-f(\cdot,\iota_0)\Big\rangle d\big\langle M^\varsigma_{\iota_0j_0}\big\rangle(s)\notag\\
&\hspace{2cm}+\Big\langle\eta(t_1),f(\cdot,\iota_1)-f(\cdot,\iota_0)\Big\rangle\notag\\
&=\big\langle\mu_0,f(\cdot,\iota_0)\big\rangle+\int_0^{t_1}\Big\langle\eta(s),\hat{\mathcal{L}}_{\iota_0}\big(\eta(s)\big)f(\cdot,\iota_0)\Big\rangle ds+\Big\langle\eta(t_1),f(\cdot,\iota_1)-f(\cdot,\iota_0)\Big\rangle,\label{eq-SMV-varsig-I2}
\end{align}
which implies that \eqref{eq-SMV-varsig-I1} also holds  for $t=t_1$. In view of \thmref{thm-Oelschlager},
 for any $0\le s\le t_1$, $\eta(s)=\mathscr{L}(z_1(s))$, which is the unique solution to the following stochastic differential equation
\begin{align*}
dz_1(s)&=b\Big(z_1(s),\mathscr L\big(z_1(s)\big),\iota_0\Big)ds+\sigma\Big(z_1(s),\mathscr L\big(z_1(s)\big),\iota_0\Big)d\tilde w(s),\quad 0\le s\le t_1,\\
\mathscr L\big(z_1(0)\big)&=\mu_0.
\end{align*}
Likewise, we can show that $\eta(t)$ satisfies \eqref{eq-SMV-varsig-k} for any $t_{k-1}\le t\le t_k$. Hence, according to \thmref{thm-Oelschlager} again, $\eta(s)=\mathscr{L}(z_k(s))$ which is the unique solution to the following equation
\begin{align*}
dz_k(s)&=b\Big(z_k(s),\mathscr L\big(z_k(s)\big),\iota_{k-1}\Big)ds+\sigma\Big(z_k(s),\mathscr L\big(z_k(s)\big),\iota_{k-1}\Big)d\tilde w(s),\quad t_{k-1}\le s\le t_k,\\
\mathscr L\big(z_k(t_{k-1})\big)&=\mathscr L\big(z_{k-1}(t_{k-1})\big).
\end{align*}
As a consequence, for $0\le s\le r$,  $\eta(s)=\mathscr{L}(z(s))$ is the unique solution to the equation
\begin{align*}
dz(s)&=b\Big(z(s),\mathscr L\big(z(s)\big),\varsigma(s_-)\Big)ds+\sigma\Big(z(s),\mathscr L\big(z(s)\big),\varsigma(s_-)\Big)d\tilde w(s),\quad 0\le s\le r,\\
\mathscr L\big(z(0)\big)&=\mu_0.
\end{align*}
It is clear that $\eta(\cdot)\in C\big([0,r],\mathscr{M}_1\big)$. This completes the Step 1.

Since for each $r>0$, $\mu_N(r)$ depends on the history of the switching process $\alpha(t)$ for $t\in[0,r)$, it is more convenient to consider equation \eqref{eq-SMV-varsig} for $0\le t<r$. Similar to Step 1, that we can define a mapping $\Lambda_r:D_f\big([0,r),\mathbb{S}\big)\to D\big([0,r),\mathscr{M}_1\big)$ that maps each $\varsigma(\cdot)\in D_f\big([0,r),\mathbb{S}\big)$ to the unique solution $\eta(\cdot)\in D\big([0,r),\mathscr{M}_1\big)$ to  equation \eqref{eq-SMV-varsig} for $0\le t<r$.
It also follows from Step 1, \eqref{eq-SMV-varsig-I1}, and \eqref{eq-SMV-varsig-I2} that the  unique solution $\eta(\cdot)\in D\big([0,r],\mathscr{M}_1\big)$ to  equation \eqref{eq-SMV-varsig} for $0\le t\le r$ satisfies
\begin{equation}\label{eqn-endpoint of SMV}
\eta(r)=\Lambda_r\varsigma(r_-)=\lim_{s\to r_-}\Lambda_r\varsigma(s),
\end{equation} for any value of $\varsigma(r)$.

For each $0<r_1\le r_2$, denote the truncation mappings $\Pi_{r_2,r_1}^{\mathbb{S}}: D\big([0,r_2),\mathbb{S}\big)\to D\big([0,r_1),\mathbb{S}\big)$ and $\Pi_{r_2,r_1}^{\mathscr{M}_1}: D\big([0,r_2),\mathscr{M}_1\big)\to D\big([0,r_1),\mathscr{M}_1\big)$ by
$$
\Pi_{r_2,r_1}^{\mathbb{S}}\varsigma(s)=\varsigma(s),\quad \Pi_{r_2,r_1}^{\mathscr{M}_1}\eta(s)=\eta(s),
$$
for all $\varsigma\in D\big([0,r_2),\mathbb{S}\big)$, $\eta\in D\big([0,r_2),\mathscr{M}_1\big)$   and  $0\le s<r_1$.
Then we have the following lemma that shows the continuity and consistency of the mapping $\Lambda_r$. To keep the continuity of the flow of presentation,
we relegate the proof to the Appendix.

\begin{lem}\label{lem-Lambda}
The following assertions hold.
\begin{itemize}
  \item[{\rm(i)}] The mapping $\Lambda_r:D_f\big([0,r),\mathbb{S}\big)\to C\big([0,r),\mathscr{M}_1\big)$ is continuous for any $0<r\le T$.
  \item[{\rm(ii)}] For any $0<r\le T$ and $\varsigma\in D_f\big([0,r_2),\mathbb{S}\big)$, the following consistent identity holds
\begin{equation}\label{eqn-consistency}
\Lambda_{r_1}\Pi_{r_2,r_1}^{\mathbb{S}}\varsigma=\Pi_{r_2,r_1}^{\mathscr{M}_1}\Lambda_{r_2}\varsigma.
\end{equation}
\end{itemize}
\end{lem}

Next, for $r>0$, we define $\alpha_{r_-}:\Omega\to D_f\big([0,r),\mathbb{S}\big)$ by $\alpha_{r_-}(s,\omega)=\alpha(s,\omega)$ for $0\le s<r$ and $\omega\in\Omega$.  Let $y_r(\cdot)$ be the solution to the following equation
$$
dy_r(s)=b\Big(y_r(s),\Lambda_r(\alpha_{r_-})(s),\alpha_{r_-}(s_-)\Big)ds+\sigma\Big(y_r(s),\Lambda_r(\alpha_{r_-})(s),\alpha_{r_-}(s_-)\Big)d\tilde w(s),\quad 0\le s< r,
$$
with $y_r(0)=y(0)$ such that $\mathscr{L}(y(0))=\mu_0$. It follows from \lemref{lem-Lambda} that $y_{r_1}(s)=y_{r_2}(s)$ for $0<s<r_1<r_2$.  Hence, we can define $y(s)=y_r(s)$ for $0\le s<r$ and obtain
\begin{equation}\label{eqn-y-Lambda}
dy(s)=b\Big(y(s),\Lambda_r(\alpha_{r_-})(s),\alpha_{r_-}(s_-)\Big)ds+\sigma\Big(y(s),\Lambda_r(\alpha_{r_-})(s),\alpha_{r_-}(s_-)\Big)d\tilde w(s),\quad 0\le s< r.
\end{equation}

\underline{Step 2:}
Prove that  $\Lambda_r(\alpha_{r_-})(r_-)=\mathscr{L}\big(y(r)\big|\mathcal{F}^\alpha_{r_-}\big)$ for each $r>0$.

First, we show $\Lambda_r(\alpha_{r_-})(r_-)=\mathscr{L}\big(y(r_-)\big|\mathcal{F}^\alpha_{r_-}\big)$ for each $0<r\le T$. Note that according to Step 1, on the set $\alpha_{r_-}=\varsigma\in D_f\big([0,r),\mathbb{S}\big)$ we have $\Lambda_r(\alpha_{r_-})(s)=\mathscr{L}\big(y(s)\big)$ for any $0\le s<r$. Therefore,
$\Lambda_r(\alpha_{r_-})(r_-)=\mathscr{L}\big(y(r_-)\big).$  As a consequence, for any $f\in C_b\big(\mathbb{R}^d\big)$,
\begin{equation}\label{eqn-consistent-condExpt}
\EE\Big(\big\langle\Lambda_r(\alpha_{r_-})(r_-),f\big\rangle\Big|\sigma\big\{\alpha_{r_-}\big\}\Big)=\EE\big(f(y(r_-))\big|\sigma\big\{\alpha_{r_-}\big\}\big).
\end{equation}

Next, let $\mathscr{B}^{\mathbb{S}}_{r_-}$ and $\mathscr{B}^{\mathscr{M}_1}_{r}$ be the Borel $\sigma-$fields on $D\big([0,r),\mathbb{S}\big)$ and $D\big([0,r),\mathscr{M}_1\big)$, respectively. For each $0<s\le r \le T$ denote the mappings  $\pi_{r,s}^{\mathbb{S}}: D\big([0,r),\mathbb{S}\big)\to \mathbb{S}$ and $\pi_{r,s}^{\mathscr{M}_1}: D\big([0,r),\mathscr{M}_1\big)\to \mathscr{M}_1$ by
$$
\pi_{r,s}^{\mathbb{S}}\varsigma=\varsigma(s),\quad \pi_{r,s}^{\mathscr{M}_1}\eta=\eta(s)\quad\text{for all }\varsigma\in D\big([0,r),\mathbb{S}\big), \eta\in D\big([0,r),\mathscr{M}_1\big).
$$
Then we have
$$
\mathscr{B}^{\mathbb{S}}_{r_-}=\sigma\big\{\pi_{r,s}^{\mathbb{S}}:0\le s<r\big\},\quad \mathscr{B}^{\mathscr{M}_1}_{r}=\sigma\big\{\pi_{r,s}^{\mathscr{M}_1}:0\le s\le r\big\}.
$$
Hence,
\begin{align}
\sigma\big\{\alpha_{r_-}\big\}=\alpha_{r_-}^{-1}\big(\mathscr{B}^{\mathbb{S}}_{r_-}\big)&=\alpha_{r_-}^{-1}\big( \sigma\big\{\pi_{r,s}^{\mathbb{S}}:0\le s<r\big\}\big)\notag\\
&=\sigma\big\{\big[\big(\pi_{r,s}^{\mathbb{S}}\big)\circ\alpha_{r_-}\big]^{-1} :0\le s<r\big\}=\mathcal{F}^\alpha_{r_-}.\label{eqn-consistent-sigma}
\end{align}

Since $\Lambda_r:D_f\big([0,r),\mathbb{S}\big)\to C\big([0,r),\mathscr{M}_1\big)$ is continuous, $\Lambda_r(\alpha_{r_-})(r_-)$ is $\sigma\big\{\alpha_{r_-}\big\}$-measurable.  In addition, $\big(\mathscr{M}_1,\|\cdot\|_{BL}\big)$ is equivalent to the space $\mathcal{P}(\mathbb{R}^d)$ equipped with the weak topology. Therefore, it follows from \eqref{eqn-consistent-condExpt} and \eqref{eqn-consistent-sigma} that
$$
\big\langle\Lambda_r(\alpha_{r_-})(r_-),f\big\rangle=\EE\big(f(y(r_-))\big|\mathcal{F}^\alpha_{r_-}\big) \quad \text{for all } f\in C_b\big(\mathbb{R}^d\big).
$$
This implies $\Lambda_r(\alpha_{r_-})(r_-)=\mathscr{L}\big(y(r_-)\big|\mathcal{F}^\alpha_{r_-}\big)$. By part (ii) of \lemref{lem-Lambda}, $\Lambda_r(\alpha_{r_-})(s)=\Lambda_T(\alpha_{T_-})(s)$ for any $0<s<r\le T$. Thus, by \lemref{lem-Lambda}(i), for each $0<s<T$, we obtain
$$
\Lambda_T(\alpha_{T_-})(s)=\Lambda_T(\alpha_{T_-})(s_-)=\Lambda_s(\alpha_{s_-})(s_-)=\mathscr{L}\big(y(s_-)\big|\mathcal{F}^\alpha_{s_-}\big).
$$
Taking $r=T$ in \eqref{eqn-y-Lambda} we have obtain
\begin{equation}\label{eqn-y-MV-minus}
dy(s)=b\Big(y(s),\mathscr{L}\big(y(s_-)\big|\mathcal{F}^\alpha_{s_-}\big),\alpha(s_-)\Big)ds+\sigma\Big(y(s),\mathscr{L}\big(y(s_-)\big|\mathcal{F}^\alpha_{s_-}\big),\alpha(s_-)\Big)d\tilde w(s),\quad 0\le s<T.
\end{equation}
Similar to \lemref{lem-Lambda}(i), we can show that solution $y(s)$ of \eqref{eqn-y-MV-minus} is continuous. This gives $\mathscr{L}\big(y(s)\big|\mathcal{F}^\alpha_{s_-}\big)=\mathscr{L}\big(y(s_-)\big|\mathcal{F}^\alpha_{s_-}\big)$ and $y(s)$ satisfies
\begin{equation}\label{eqn-y-MV}
dy(s)=b\Big(y(s),\mathscr{L}\big(y(s)\big|\mathcal{F}^\alpha_{s_-}\big),\alpha(s_-)\Big)ds+\sigma\Big(y(s),\mathscr{L}\big(y(s)\big|\mathcal{F}^\alpha_{s_-}\big),\alpha(s_-)\Big)d\tilde w(s),\quad 0\le s\le T.
\end{equation}
The assertion of Step 2 is therefore complete.

It follows from Step 1 and Step 2 (with $\varsigma$ is replaced by sample path $\alpha_{T_-}$) that the solution $\mu$ to \eqref{eq-perturb-SMV} satisfies
$$
\mu(s)=\Lambda_T(\alpha_{T_-})(s)= \mathscr{L}\big(y(s)\big|\mathcal{F}^\alpha_{s_-}\big), \quad 0<s<T.
$$
In view of \eqref{eqn-endpoint of SMV} and Step 2,
$$
\mu(T)=\Lambda_T(\alpha_{T_-})(T_-)=\mathscr{L}\big(y(T)\big|\mathcal{F}^\alpha_{T_-}\big).
$$
These imply $\mu(s)= \mathscr{L}\big(y(s)\big|\mathcal{F}^\alpha_{s_-}\big)$ for all $0\le s\le T$.

\underline{Step 3:} Prove the uniqueness of the solution of \eqref{eqn-y-MV}.

Suppose that $y_1,y_2$ are two solutions to  equation \eqref{eqn-y-MV} with same initial value $y(0)$. Denote $\mu_i(s)=\mathscr{L}\big(y_i(s)\big|\mathcal{F}^\alpha_{s_-}\big)$ for $0\le s\le T$, $i=1,2$. Then
$$
y_i(t)=y(0)+\int_0^tb\big(y_i(s),\mu_i(s),\alpha(s_-)\big)ds+\int_0^t\sigma\big(y_i(s),\mu_i(s),\alpha(s_-)\big)d\tilde w(s),\quad 0\le t\le T.
$$
Similar to \eqref{inq-|b1-b2|}, we have
\begin{align}
&\Big|b\big(y_1(s),\mu_1(s),\alpha(s_-)\big)-b\big(y_2(s),\mu_2(s),\alpha(s_-)\big)\Big|^2\notag\\
&\hspace{4cm}+\Big|\sigma\big(y_1(s),\mu_1(s),\alpha(s_-)\big)-\sigma\big(y_2(s),\mu_2(s),\alpha(s_-)\big)\Big|^2\notag\\
& \quad \le 2\Big[\big|y_1(s)-y_2(s)\big|^2+\EE\Big(\big|y_1(s)-y_2(s)\big|^2\Big|\mathcal{F}^\alpha_{s_-}\Big)\Big].\label{inq-|b1-b2|-cond}
\end{align}
Thus, by Cauchy-Schwarz and Bulkhoder-Davis-Gundy inequalities, we have
\begin{align}
\EE\big|y_1(t)-y_2(t)\big|^2&
\le2T\int_0^t\EE\Big|b\big(y_1(s),\mu_1(s),\alpha(s_-)\big)-b\big(y_2(s),\mu_2(s),\alpha(s_-)\big)\Big|^2ds\notag\\
&\quad +2\EE\Bigg|\int_0^t\Big[\sigma\big(y_1(s),\mu_1(s),\alpha(s_-)\big)-\sigma\big(y_2(s),\mu_2(s),\alpha(s_-)\big)\Big]d\tilde w(s)\Big|^2\notag\\
&\le 2T\int_0^t\EE\bigg[\big|y_1(s)-y_2(s)\big|^2+\EE\Big(\big|y_1(s)-y_2(s)\big|^2\Big|\mathcal{F}^\alpha_{s_-}\Big)\bigg]ds\notag\\
&\quad +C\int_0^t\EE\Big|\sigma\big(y_1(s),\mu_1(s),\alpha(s_-)\big)-\sigma\big(y_2(s),\mu_2(s),\alpha(s_-)\big)\Big|^2ds\notag\\
&\le C\int_0^t\EE\bigg[\big|y_1(s)-y_2(s)\big|^2+\EE\Big(\big|y_1(s)-y_2(s)\big|^2\Big|\mathcal{F}^\alpha_{s_-}\Big)\bigg]ds\notag\\
&=C\int_0^t\EE\big|y_1(s)-y_2(s)\big|^2ds.
\end{align}
In view of the Gronwall inequality, $\EE\big|y_1(t)-y_2(t)\big|^2=0$ for $0\le t\le T$. This implies that $y_1(t)=y_2(t)$ a.s.
\qed

\para{Proof of \thmref{thm-main1}}. Since $\mathscr{P}_N$ is the distribution of $\big(\mu_N(\cdot),\alpha(\cdot)\big)$, $\mathscr{P}\big(\mathcal{C}\big)=\PP\big(\alpha_T\in \mathcal{C}\big)$ for any measurable set $\mathcal{C}\subset D\big([0,T],\mathbb{S}\big)$ where $\alpha_{T}:\Omega\to D_f\big([0,T],\mathbb{S}\big)$ is defined by $\alpha_{T}(s,\omega)=\alpha(s,\omega)$ for $0\le s\le T$ and $\omega\in\Omega$. The convergence $\mathscr{L}\big(\mu_N(0)\big)\Rightarrow \delta_{\mu_0}$ in  $\mathcal{P}\big(\mathscr{M}_1,\|\cdot\|_{BL}\big)$ implies that $\mathscr{P}\big(\eta(0)=\mu_0\big)=1$.

Let $\bar \Lambda_T:D_f\big([0,T],\mathbb{S}\big)\to D\big([0,T],\mathscr{M}_1\big)$ be the mapping that maps each $\varsigma(\cdot)\in D_f\big([0,T],\mathbb{S}\big)$ to the unique solution $\eta(\cdot)\in D\big([0,T],\mathscr{M}_1\big)$ to  equation \eqref{eq-SMV-varsig} (or equivalently,
\eqref{eq-stochastic-MV}) for $0\le t\le T$. Similar to \lemref{lem-Lambda} we can show that $\bar \Lambda_T:D_f\big([0,T],\mathbb{S}\big)\to C\big([0,T],\mathscr{M}_1\big)$ is continuous.

Denote $\Gamma:D_f\big([0,T],\mathbb{S}\big)\to D\big([0,T],\mathscr{M}_1\times\mathbb{S}\big)$ by $\Gamma\varsigma=\big(\bar \Lambda_T\varsigma,\varsigma\big)$ and let $\mathcal{S}$ be the set of all pairs $(\eta,\varsigma)\in  C\big([0,T],\mathscr{M}_1\big)\times D_f\big([0,T],\mathbb{S}\big)$ satisfying equation \eqref{eq-stochastic-MV}. Since  $C\big([0,T],\mathscr{M}_1\big)\times D_f\big([0,T],\mathbb{S}\big)$ is a closed subset of $D\big([0,T],\mathscr{M}_1\times\mathbb{S}\big)$,
$$
\mathscr{P}\Big(C\big([0,T],\mathscr{M}_1\big)\times D_f\big([0,T],\mathbb{S}\big)\Big)\ge\limsup_{N\to\infty}\mathscr{P}_N\Big(C\big([0,T],\mathscr{M}_1\big)\times D_f\big([0,T],\mathbb{S}\big)\Big)=1.
$$
Thus, $\mathscr{P}\big(C\big([0,T],\mathscr{M}_1\big)\times D_f\big([0,T],\mathbb{S}\big)\big)=1$. This together with \thmref{thm-limit-characterization} implies that $\mathscr{P}(\mathcal{S})=1$. According to Step 1 in the proof of \thmref{thm-SMV}, $\mathcal{S}=\big\{\Gamma\varsigma:\varsigma \in D_f\big([0,T],\mathbb{S}\big)\big\}$. Therefore, for each measurable set $\mathcal{A}\subset D\big([0,T],\mathscr{M}_1\times\mathbb{S}\big)$ we have
$$
\mathscr{P}(\mathcal{A})=\mathscr{P}\big(\mathcal{A}\cap \mathcal{S}\big)=\mathscr{P}\Big(\varsigma\in \Gamma^{-1}(\mathcal{A})\Big)=\PP\Big(\alpha_T\in\Gamma^{-1}(\mathcal{A})\Big)=\PP\Big(\big(\bar\Lambda_T\alpha_T,\alpha_T\big)\in \mathcal{A}\Big).
$$
Since $\bar\Lambda_T\alpha_T(s)=\mathscr{L}\big(y(s)\big|
\mathcal{F}^\alpha_{s_-}\big)$, where $y(t)$ is the unique solution to \eqref{eqn-y-MV}, the above identities imply that $\mathscr{P}$ is the distribution of the process $\big(\mathscr{L}\big(y(s)\big|
\mathcal{F}^\alpha_{s_-}\big),\alpha(s)\big)$ on $D\big([0,T],\mathscr{M}_1\times\mathbb{S}\big)$. This completes the proof.
\qed

\section{$N$-Particle Mean-Field Models with
Two-Time-Scale
Markovian Switching Process}\label{sec:two-time}

As alluded to in the introduction, this section is devoted to the case that the number of particles $N\to \infty$, meanwhile, the Markov chain displays weak and strong interactions reflected by use of a small parameter $\epsilon\to 0$. We require that
$N\wedge (1/\epsilon)\to \infty$ as $\epsilon \to 0$ and $N\to \infty$.

\subsection{Formulation}
We consider a class of mean-field processes, in
which the random
switching process changes much faster than the
continuous state (or the switching process jump change much more
frequently). The basic idea is that there are inherent two-time
scales. Our interest focuses on the limit behavior of the resulting
process. Suppose that $\epsilon > 0$ is a small parameter and the
system of mean-field equations is given by
\begin{equation}\label{Eqn-Dyn-2TimeScale}
    dx^\epsilon_i(t) = b\bigg(x^\epsilon_i(t),{1\over N}\sum_{j=1}^N\delta_{x^\epsilon_j(t)},\alpha^\epsilon(t_-)\bigg)dt + \sigma\bigg(x^\epsilon_i(t),{1\over N}\sum_{j=1}^N\delta_{x^\epsilon_j(t)},\alpha^\epsilon(t_-)\bigg) dw_i(t),
\end{equation}
where $w_1(\cdot)$, $w_2(\cdot)$, $\ldots$, $w_N(\cdot)$ are independent $d$-dimensional
standard Brownian motions, and $\alpha^\epsilon(t)$ is a Markov chain with state space $\mathbb{S} =
\big\{1, \ldots, m_0\big\}$ satisfying
\begin{equation}
\PP\Big( \alpha^\epsilon(t + \Delta t) = j_0 \Big| \alpha^\epsilon (t) = i_0; x^\epsilon_{1}(s), \ldots, x^\epsilon_{N}(s), 0\leq s \leq t \Big) = q^\epsilon_{i_0j_0}\Delta t + o\big(\Delta t\big), \ i_0 \neq j_0,\end{equation}
as $\Delta t\to 0$, where the generator
$$
Q^\epsilon=\Big(q^\epsilon_{i_0j_0}\Big)_{i_0,j_0\in\mathbb{S}}={1\over\epsilon}\tilde Q+\hat Q,
$$
satisfies $q^\epsilon_{i_0j_0}\ge0$ for $i_0 \neq j_0$ and $\sum_{j_0
\in \mathcal{M}} q^\epsilon_{i_0j_0} = 0$ for each $i_0 \in \mathbb{S}$.

The model above is motivated by the work of two-time-scale Markov chains \cite{YinZhang13}. Such two-time scale Markov chains have been used widely, especially in networked systems; see for example, the manufacturing systems given in \cite{SethiZ}. It is readily seen that the Markov chain has a fast varying part and a slowly changing part.
Suppose
$$
\tilde Q=\text{diag}\Big[\tilde Q^1,\tilde Q^2,\ldots,\tilde Q^l\Big].
$$
Then the state space $\mathbb S$ of the underlying Markov chain is decomposable into $l$ subspaces.
These subspaces are not completely separated. There are weak interactions among the subspaces due to the use of the slowly varying part of the generator $\hat Q$. Such a structure is often referred to as nearly decomposable Markov chain.

We relabel the states so that
$$
\mathbb S=\mathbb S_1\cup\mathbb S_2\cup\ldots\cup\mathbb S_l,
$$
with $\mathbb S_i=\big\{s_{i1},s_{i2},\ldots,s_{im_i}\big\}$ and $m_0=m_1+m_2+\ldots+m_l$ such that $\tilde Q^i$, the generator associated with the subspace $\mathbb S_i$ for each $i=1,\ldots,l$. Assume that each $\tilde Q^i$ is irreducible. As a consequence, the corresponding $\mathbb S_i$ for $i=1,\ldots,l$ consist of recurrent states belonging to $l$ ergodic classes. Let $\nu^i=\big(\nu_{s_{i1}},\nu_{s_{i2}},\ldots,\nu_{s_{im_i}}\big)$ be the stationary distribution corresponding to $\tilde Q^i$, $1\le i\le l$, and $\tilde\nu=\text{diag}\big[\nu^1,\nu^2,\ldots,\nu^l\big]\in\mathbb{R}^{l\times m_0}$.
Following the ideas in \cite{YinZhang13}, we aim to reduce the computational complexity. The rationale is that we take advantage of the fast and slow motions and strong and weak interactions of the systems so that we can naturally divide the state space of the switching process into subsystems or groups. Within each subsystem, the states look alike in that they vary at the same speed, and among different subsystems, the variations are relatively slowly.
To proceed, we lump the states of the jump component in each $\mathbb S_i$ into a single state and define
$$
\bar\alpha^\epsilon(t)=i\quad\text{if}\quad\alpha^\epsilon(t)\in\mathbb{S}_i.
$$
Denote the state space of $\bar\alpha^\epsilon(\cdot)$ by $\bar {\mathbb S}=\big\{1,2,\ldots,l\big\}$.   It follows from \cite[Theorem 5.27]{YinZhang13} that $\bar\alpha^\epsilon(\cdot)$ converges weakly to $\bar\alpha(\cdot)$, a Markov chain with the state space $\bar{\mathbb{S}}$ and generator  $\bar Q$ defined by
\begin{equation}\label{Def-bar-Q}
\bar Q=\big(\bar q_{ij}\big)_{l\times l}=\tilde\nu\hat Q {1\!\!1},
\end{equation}
where $1\!\!1=\text{diag}\big[{1\!\!1}_{m_1},{1\!\!1}_{m_2},\ldots,{1\!\!1}_{m_l}\big]\in\mathbb{R}^{m_0\times l}$ and ${1\!\!1}_k=\big(1,1,\ldots,1\big)'\in{\mathbb R}^k$.

For $(x,\mu,i)\in\mathbb{R}^d\times\mathscr{M}_1\times\bar{\mathbb{S}}$, denote $a\big(x,\mu,s_{ij}\big)=\sigma\big(x,\mu,s_{ij}\big)\sigma'\big(x,\mu,s_{ij}\big)$ and
\begin{equation}\label{Def-bar-ab}
\bar b\big(x,\mu,i\big)=\sum_{j=1}^{m_{i}}\nu_{s_{ij}}b\big(x,\mu,s_{ij}\big),\quad\bar a\big(x,\mu,i\big)=\sum_{j=1}^{m_{i}}\nu_{s_{ij}}a\big(x,\mu,s_{ij}\big).
\end{equation}
For simplicity, we assume that the initial values $x^\epsilon_i(0)=x_{0,i}$, $i=1,2,\ldots$, are independent of $\epsilon$ and that $x_{0,i}$ are independent of $\alpha^\epsilon(\cdot)$. We make the following assumption.
\begin{itemize}
\item[(A3)] The matrix-valued function $\bar a(\cdot,\cdot,\cdot)$ on ${\mathbb{R}}^d\times \mathscr{M}_1\times \bar{\mathbb{S}} \to \mathbb{R}^{d\times d}$ has a representation $\bar a(x,\mu,i)=\bar \sigma(x,\mu,i)\bar \sigma'(x,\mu,i)$ where $\bar \sigma(\cdot,\cdot,\cdot)$ is bounded and there is a constant $L$ such that
    $$
    \big|\bar \sigma(x,\mu,i)-\bar \sigma(y,\eta,i)\big|\le L\big(|x-y|+\|\mu-\eta\|_{BL}\big),
    $$
    for all $x,y\in\mathbb{R}^d$, $\mu,\eta\in\mathscr{M}_1$  and $i\in\bar{\mathbb{S}}$.
\end{itemize}

For $N\ge1$, $\epsilon>0$, $0\le t\le T$, and $A\in\mathcal{B}(\mathbb{R}^d)$ denote
\begin{equation}
    \mu^\epsilon_{N}(t,A) = \frac{1}{N}\sum_{j = 1}^N \delta_{x^\epsilon_j(t)}(A).
\end{equation}
Then $\mu^\epsilon_N(t)$, $0\le t\le T$, defines a process on the space  $\mathscr{M}_1$ of probability measures on $\mathbb{R}^d$. Because of the assumption on the initial values of $x^\epsilon_i(0)$, $\mu^\epsilon_N(0)=\mu_N(0)$ does not depend on $\epsilon$. Let $d_{\bar{\mathbb{S}}}$ and $\bar d$, respectively, be the metrics on $\bar{\mathbb{S}}$ and $\mathscr{M}_1\times\bar{\mathbb{S}}$ defined by similar way to $d_{\mathbb{S}}$ and $d$ in \eqref{Def-metric-d}. Denote $\bar{\mathscr{P}}_N^\epsilon$ by the induced probability measure of $\big(\mu_N^\epsilon(\cdot),\bar{\alpha}^\epsilon(\cdot)\big)$ on $D\big([0,T],\mathscr{M}_1\times\bar{\mathbb{S}}\big)$. We can show that $\big(\mu_N^\epsilon(\cdot),\bar{\alpha}^\epsilon(\cdot)\big)\in C([0,T],\mathscr M_1)\times D_f([0,T],\bar{\mathbb{S}})$, a closed subspace of $D\big([0,T],\mathscr{M}_1\times\bar{\mathbb{S}}\big)$.

\begin{thm}\label{thm-main2}
Assume {\rm(A1)}, {\rm(A2)}, {\rm(A3)}, and
$$
\sup_{N\in\mathbb N}\EE\big\langle \mu_N(0),\psi\big\rangle<\infty,\ \mathscr{L}\big(\mu_N(0)\big)\Rightarrow \delta_{\mu_0}\text{ in }\mathcal{P}\big(\mathscr{M}_1,\|\cdot\|_{BL}\big).
$$
Then $\big(\mu_N^\epsilon(\cdot),\bar\alpha^\epsilon(\cdot)\big)$ converges weakly to process $\big(\bar \mu_{\bar\alpha}(\cdot),\bar\alpha(\cdot)\big)$ as $\epsilon\to0$ and $N\to\infty$ satisfying
$(1/\epsilon)\wedge N\to\infty$, where
$$
\big(\bar\mu_{\bar\alpha}(t),\bar\alpha(t)\big)=\Big(\mathscr L\big(\bar y(t)\big|\mathcal{F}^{\bar\alpha}_{t_-}\big),\bar\alpha(t)\Big),\quad 0\le t\le T,
$$
and $\bar\zeta(t)$, $0\le t\le T$, is the unique solution of the following  stochastic differential equation
$$\begin{array}{ll}&\disp
d\bar y(t)=\bar b\bigg(\bar y(t),\mathscr L\big(\bar y(t)\big|\mathcal{F}^{\bar\alpha}_{t_-}\big),\bar\alpha(t_-)\bigg)dt+\bar\sigma\bigg(\bar y(t),\mathscr L\big(\bar y(t)\big|\mathcal{F}^{\bar\alpha}_{t_-}),
\bar\alpha(t_-)\bigg)d\tilde w(t), \\
&\disp \mathscr L\big(\bar y(0)\big)=\mu_0, \end{array}
$$
where $\tilde w(\cdot)$ is a standard Brownian motion independent of $\bar\alpha(\cdot)$.
\end{thm}

\subsection{Weak Compactness and Auxiliary Estimates}

For $N\ge1$, $\epsilon>0$ and $t>0$ denote
$$
\mathcal{F}^{\bar \alpha}_{t_-}=\sigma\big\{\bar\alpha(s):0\le s<t\big\},\quad \mathcal{F}^{N,\epsilon}_{t}=\sigma\big\{w_i(s),\alpha^\epsilon(s):0\le s\le t,1\le i\le N\big\}.
$$

\begin{rem}\label{Rem-Boundedness}{\rm
It follows from the proof of Lemma \ref{lem-moments} that for each  $p$, $0<p\le1$, there exists a constant $C$ independent of $N$ and $\epsilon$ such that
\begin{equation}
\sup_{0\le t\le T}\EE\Big(\big\langle\mu_N^\epsilon(t),\psi\rangle+1\Big)^p\le C, \label{moment-inq-timescale}
\end{equation}
and, for $0\le s\le t\le T$,
\begin{equation}
e^{-C(t-s)}\Big(\big\langle\mu_N^\epsilon(s),\psi\rangle+1\Big)^p\le \EE\Big[\Big(\big\langle\mu_N^\epsilon(t),\psi\rangle+1\Big)^p\Big|\mathcal{F}^{N,\alpha^\epsilon}_s\Big]\le
e^{C(t-s)}\Big(\big\langle\mu_N^\epsilon(s),\psi\rangle+1\Big)^p.\label{3moment-inq-timescale}
\end{equation}
}
\end{rem}

\begin{prop}\label{Prop-weakcompact}
Assume that all the assumptions of \thmref{thm-main2} hold. Then the sequence $\{\mu^\epsilon_N(\cdot),\bar\alpha^\epsilon(\cdot)\}$ is weakly compact in the topology of weak convergence of probability measure on $D([0,T],\mathscr{M}_1\times\bar{\mathbb{S}})$.
\end{prop}

\para{Proof.}
According to inequalities \eqref{moment-inq-timescale}-\eqref{3moment-inq-timescale} and the arguments in \lemref{lem-compact-containment} and \lemref{lem-weak-compactness}, for any $\delta>0$ there exists a compact set $K_\delta$ in $(\mathscr{M}_1,\|\cdot\|_{BL})$ such that
$$
\inf_{\epsilon;N}\PP\Big(\mu_N^\epsilon(t)\in K^\delta_\delta\,\,\,\forall \,\,0\le t\le T\Big)\ge1-\delta,
$$
and there exists a constant $C$ independent of $N$ and $\epsilon$ such that
$$
\EE\Big[\big\|\mu_N^\epsilon(t+s)-\mu_N^\epsilon(t)\big\|_{BL}\Big|\mathcal{F}^{N,\alpha^\epsilon}_t\Big]\le C\sqrt{\delta} \EE\Big[\big\langle\mu_N^\epsilon(T),\psi\big\rangle+1\Big|\mathcal{F}^{N,\alpha^\epsilon}_t\Big].
$$
Since $\bar \alpha^\epsilon(\cdot)$ converges weakly to $\bar\alpha(\cdot)$ (see Theorem 7.4 \cite{YinZhang13}) we obtain the compactness of $(\mu_N^\epsilon(\cdot),\bar\alpha^\epsilon(\cdot))$.
\qed

For each function $f(\cdot,\cdot)$ with $f(\cdot,s_{ij})\in C^2(\mathbb{R}^d)$ for $s_{ij}\in\mathbb{S}$, denote the operator associated to \eqref{Eqn-Dyn-2TimeScale} by
\begin{equation}\label{Def-Lf-timescale}
\mathcal{L}^\epsilon(\mu)f\big(x,s_{ij}\big)=b'\big(x,\mu,s_{ij}\big)\nabla_{x}f\big(x,s_{ij}\big)+{1\over2}\Big(a\big(x,\mu,s_{ij}\big)\nabla_x\Big)'\nabla_xf\big(x,s_{ij}\big)+Q^\epsilon f(x)(s_{ij}),
\end{equation}
where
$$Q^\epsilon f(x)(s_{ij})=\sum_{s_{i_1j_1}\in \mathbb{S}}q^\epsilon_{s_{ij},s_{i_1j_1}}\big(f\big(x,s_{i_1j_1}\big)-f\big(x,s_{ij}\big)\big).$$
To approximate this operator for small values of $\epsilon$ we define for each function $g(\cdot,\cdot)$ such that $g(\cdot,i)\in C^2(\mathbb{R}^d)$ for $i\in\bar{\mathbb{S}}$,
\begin{equation}\label{Def-bar-L-timescale}
\bar{\mathcal{L}}(\mu)g\big(x,i\big)=\bar b'\big(x,\mu,i\big)\nabla_{x}g\big(x,i\big)+{1\over2}\Big(\bar a(x,\mu,i)\nabla_x\Big)'\nabla_xg\big(x,i\big)+\bar Qg(x)(i),
\end{equation}
 where $\bar Qg(x)(i)=\sum_{j\in\bar{\mathbb{S}}}\bar q_{ij}\big(g(x,j)-g(x,i)\big)$.

Note that for $f(\cdot)\in C^2(\mathbb{R}^d)$, we can define $\mathcal{L}^\epsilon(\mu)f\big(x,s_{ij}\big)$ as in \eqref{Def-Lf-timescale}  by taking $f(x,s_{ij})\equiv f(x)$ for $x\in\mathbb{R}^d$ and $s_{ij}\in\mathbb{S}$. Similarly, for each $f(\cdot)\in C^2(\mathbb{R}^d)$, we can define
 $\bar{\mathcal{L}}(\mu)f\big(x,i\big)=\bar{\mathcal{L}}(\mu)g\big(x,i\big)$
  where $g(x,i)\equiv f(x)$ for any $x\in\mathbb{R}^d$
  and $i\in\bar{\mathbb{S}}$. We have the following approximation. To make the presentation more transparent, its proof is given in the Appendix.

\begin{lem}\label{Lem-operator-est} Under Assumption {\rm A}, for any $f\in C^3_c(\mathbb{R}^d)$ there is a constant $C$ independent of $N$ and $\epsilon$ such that
\begin{equation}
\sup_{0\le t\le T} \EE\bigg|\int_0^t\Big\langle\mu_N^\epsilon(s),\mathcal{L}^\epsilon\big(\mu_N^\epsilon(s)\big)f\big(\alpha^\epsilon(s)\big)\Big\rangle ds-\int_0^t\Big\langle\mu_N^\epsilon(s),\bar{\mathcal{L}}\big(\mu_N^\epsilon(s)\big)f\big(\bar\alpha^\epsilon(s)\big)\Big\rangle ds\bigg|\le C\epsilon^{1/6}.\label{operator-est}
\end{equation}
\end{lem}

\subsection{Weak Convergence and Stochastic McKean-Vlasov Equation with Two-Time-Scale Markovian Switching}
Similar to \eqref{Def-Mart of MC}, we define the martingale  associated with  the limiting Markovian switching process $\bar\alpha(\cdot)$ by
\begin{equation}\label{Def-Mart of MC bar}
\bar M_{ij}(t)=\big[\bar M_{ij}\big](t)- \big\langle \bar M_{ij}\big\rangle(t),\quad i,j\in\bar{\mathbb{S}},\, 0\le t\le T,
\end{equation}
where
$
\big[\bar M_{ij}\big](t)=\sum_{0\le s\le t}{1\!\!1}\big(\bar\alpha(s_-)=i\big){1\!\!1}\big(\bar\alpha(s)=j\big)$  and $\big\langle\bar M_{ij}\big\rangle(t) =\int_0^t\bar q_{i_0j_0}{1\!\!1}\big(\bar\alpha(s_-)=i_0\big)ds$.
In addition, by a similar way to \eqref{Def-sample Mart of MC}, we define sample path of the martingale  associated with a sample path $\bar \varsigma \in D_f\big([0,T],\bar{\mathbb{S}}\big)$ of  $\bar\alpha(t)$ by
\begin{equation}\label{Def-bar-sample Mart of MC}
\bar M_{ij}^{\bar \varsigma}(t)=\big[\bar M_{ij}^{\bar \varsigma}\big](t)- \big\langle\bar  M_{ij}^{\bar \varsigma}\big\rangle(t),\qquad i\ne j\in \bar{\mathbb{S}},\,0\le t\le T,
\end{equation}
where $\big[\bar M_{ij}^{\bar\varsigma}\big](t)=\sum_{0\le s\le t}{1\!\!1}\big(\bar \varsigma(s_-)=i\big){1\!\!1}\big(\bar \varsigma(s)=j\big)$ and $\big\langle \bar M_{ij}^{\bar \varsigma}\big\rangle(t) =\int_0^t\bar q_{ij}{1\!\!1}\big(\bar \varsigma(s_-)=i\big)ds$. In addition, for $\bar\varsigma\in D\big([0,T],\bar{\mathbb{S}}\big)\backslash D_f\big([0,T],\bar{\mathbb{S}}\big)$ and $0\le t\le T$ define $\bar M_{ij}^{\bar \varsigma}(t)=\big[\bar M_{ij}^{\bar \varsigma}\big](t)=\big\langle\bar  M_{ij}^{\bar \varsigma}\big\rangle(t)=0.$

\begin{prop}\label{Prop-SDE-char} Assume {\rm (A1)}, {\rm (A2)}, and $\sup_{N\ge 1}\EE\langle\mu_N(0),\psi\rangle<\infty$.
Denote by $\bar{\mathscr{P}}$ the limit of an arbitrary weakly convergence subsequence $\bar{\mathscr{P}}_{N_k}^{\epsilon_k}$ as $k\to\infty$ where $(\epsilon_k,N_k)$
satisfies $\min\{1/\epsilon_k,N_k\}\to\infty$ as $k\to\infty$. Then for $\bar{\mathscr{P}}$-almost all $(\eta,\bar\varsigma)\in D\big([0,T],\mathscr{M}_1\times\bar{\mathbb{S}}\big)$, the equation
\begin{equation}\label{eqn-char-perb}
\begin{array}{ll}
\disp
\big\langle\eta(t),g\big(\cdot,\bar\varsigma(t)\big)\big\rangle
\!\!\disp&=\big\langle\eta(0),g\big(\cdot,\bar\varsigma(0)\big)\big\rangle
+\disp\int_0^t\Big\langle\eta(s),\bar{\mathcal{L}}\big(\eta(s)\big)g\big(\cdot,\bar\varsigma(s_-)\big)\Big\rangle ds\\
&\qquad \qquad\qquad\quad+\disp\sum_{i,j\in\bar{S}}\disp\int_0^t\Big\langle \eta(s), g(\cdot,j)-g(\cdot,i)\Big\rangle d{\bar M}^{\bar\varsigma}_{ij}(s),\end{array}
\end{equation}
holds for any test function $g(\cdot,i)\in C^2_c(\mathbb{R}^d)$, $i\in \bar{\mathbb{S}}$.
\end{prop}

\para{Proof.} Let $\bar{\mathscr{P}}$ be the weak limit on $D\big([0,T],\mathscr{M}_1\times\bar{\mathbb{S}}\big)$ of  a subsequence $\bar{\mathscr{P}}^{\epsilon_k}_{N_k}$  of $\bar{\mathscr{P}}^{\epsilon}_{N}$ as $k\to\infty$ where $(\epsilon_k,N_k)$ is satisfies $\min\{1/\epsilon_k,N_k\}\to\infty$ as $k\to\infty$.

First, we prove that for $\bar{\mathscr{P}}-$almost all $(\eta,\bar\varsigma)\in D\big([0,T],\mathscr{M}_1\times\bar{\mathbb{S}}\big)$, \eqref{eqn-char-perb} holds for $g(x,i)=f(x)$ with  $f\in C^2_c(\mathbb{R}^d)$. For each $(\eta,\bar\varsigma)$ denote
$$
\bar M_f(t)=\big\langle\eta(t),f\big\rangle-\big\langle\eta(0),f\big\rangle-\int_0^t\Big\langle\eta(s),\bar{\mathcal{L}}\big(\eta(s)\big)f\big(\cdot,\bar\varsigma(s_-)\big)\Big\rangle ds,
$$
which defines a function on $D\big([0,T],\mathscr{M}_1\times\bar{\mathbb{S}}\big)$. We observe that for a fixed pair $(\eta,\bar\varsigma)\in D\big([0,T],\mathscr{M}_1\times\bar{\mathbb{S}}\big)$, if \eqref{eqn-char-perb} holds for any test function $g(\cdot,\cdot)$ such that $g(x,i)\equiv f(x)$, $f\in C^3_c(\mathbb{R}^d)$, it also holds for any test function $g(\cdot,\cdot)$ such that $g(x,i)\equiv f(x)$, $f\in C^2_c(\mathbb{R}^d)$. Thus, in order  to prove \eqref{eqn-char-perb} holds $\bar{\mathscr{P}}-$almost surely for any $g=f\in C^2_c(\mathbb{R}^d)$, it suffices to show that
\begin{equation}\label{suffice-barP}
\bar{\mathscr{P}}\Big((\eta,\bar\varsigma):\big|\bar M_f(t)\big|=0\Big)=1\quad \forall\,\, f\in C^3_c(\mathbb{R}^d).
\end{equation}
Take $\delta>0$. Since the set $\big\{(\eta,\bar\varsigma):|\bar M_f(t)|\le\delta\big\}$ is closed in $D\big([0,T],\mathscr{M}_1\times\bar{\mathbb{S}}\big)$ we have
\begin{equation}\label{est-barP1}
\bar{\mathscr{P}}\Big(\big|\bar M_f(t)\big|\le \delta\Big)\ge\limsup_{k\to\infty}\bar{\mathscr{P}}^{\epsilon_k}_{N_k}\Big(\big|\bar M_f(t)\big|\le \delta \Big)\ge1-\liminf_{k\to\infty}\bar{\mathscr{P}}^{\epsilon_k}_{N_k}\Big(\big|\bar M_f(t)\big|>\delta\Big).
\end{equation}
Note that since $\bar{\mathscr{P}}^{\epsilon_k}_{N_k}$ is the distribution of $\big(\mu^{\epsilon_k}_{N_k}(\cdot),\bar\alpha^{\epsilon_k}(\cdot)\big)$,
\begin{align}
&\bar{\mathscr{P}}^{\epsilon_k}_{N_k}\Big(\big|\bar M_f(t)\big|>\delta\Big)\notag\\
&=\PP\bigg(\bigg|\big\langle\mu^{\epsilon_k}_{N_k}(t),f\big\rangle-\big\langle\mu^{\epsilon_k}_{N_k}(0),f\big\rangle-\int_0^t\Big\langle\mu^{\epsilon_k}_{N_k}(s),\bar{\mathcal{L}}\big(\mu^{\epsilon_k}_{N_k}(s)\big)f\big(\cdot,\bar\alpha^{\epsilon_k}(s)\big)\Big\rangle ds\bigg|>\delta\bigg)\notag\\
&\le \PP\bigg(\bigg|\big\langle\mu^{\epsilon_k}_{N_k}(t),f\big\rangle-\big\langle\mu^{\epsilon_k}_{N_k}(0),f\big\rangle-\int_0^t\Big\langle\mu^{\epsilon_k}_{N_k}(s),\mathcal{L}^{\epsilon_k}\big(\mu^{\epsilon_k}_{N_k}(s)\big)f\big(\cdot,\alpha^{\epsilon_k}(s)\big)\Big\rangle ds\bigg|>{\delta\over2}\bigg)\notag\\
&\qquad+\PP\bigg(\bigg|\int_0^t\Big\langle\mu^{\epsilon_k}_{N_k}(s),\mathcal{L}^{\epsilon_k}\big(\mu^{\epsilon_k}_{N_k}(s)\big)f\big(\cdot,\alpha^{\epsilon_k}(s)\big)\Big\rangle ds\notag\\
&\qquad\qquad\qquad\qquad\qquad\qquad\qquad-\int_0^t\Big\langle\mu^{\epsilon_k}_{N_k}(s),\bar{\mathcal{L}}\big(\mu^{\epsilon_k}_{N_k}(s)\big)f\big(\cdot,\bar\alpha^{\epsilon_k}(s)\big)\Big\rangle ds\bigg|>{\delta\over2}\bigg).\label{est-prob-timescale}
\end{align}
By Lemma \ref{Lem-operator-est}, we obtain
\begin{align}
&\PP\bigg(\bigg|\int_0^t\Big\langle\mu^{\epsilon_k}_{N_k}(s),\mathcal{L}^{\epsilon_k}\big(\mu^{\epsilon_k}_{N_k}(s)\big)f\big(\cdot,\alpha^{\epsilon_k}(s)\big)\Big\rangle ds-\int_0^t\Big\langle\mu^{\epsilon_k}_{N_k}(s),\bar{\mathcal{L}}\big(\mu^{\epsilon_k}_{N_k}(s)\big)f\big(\cdot,\bar\alpha^{\epsilon_k}(s)\big)\Big\rangle ds\bigg|>{\delta\over2}\bigg)\notag\\
&\le {2\over\delta}\EE\bigg|\int_0^t\Big\langle\mu^{\epsilon_k}_{N_k}(s),\mathcal{L}^{\epsilon_k}\big(\mu^{\epsilon_k}_{N_k}(s)\big)f\big(\cdot,\alpha^{\epsilon_k}(s)\big)\Big\rangle ds-\int_0^t\Big\langle\mu^{\epsilon_k}_{N_k}(s),\bar{\mathcal{L}}\big(\mu^{\epsilon_k}_{N_k}(s)\big)f\big(\cdot,\bar\alpha^{\epsilon_k}(s)\big)\Big\rangle ds\bigg|\notag\\
&\le {C\epsilon_k^{1/6}\over\delta}.\label{est-prob1-timescale}
\end{align}
Next, denote
$$
M^{\epsilon_k}_{N_k,f}(t)=\big\langle\mu^{\epsilon_k}_{N_k}(t),f\big\rangle-\big\langle\mu^{\epsilon_k}_{N_k}(0),f\big\rangle-\int_0^t\Big\langle\mu^{\epsilon_k}_{N_k}(s),\mathcal{L}^{\epsilon_k}\big(\mu^{\epsilon_k}_{N_k}(s)\big)f\big(\cdot,\alpha^{\epsilon_k}(s)\big)\Big\rangle ds.
$$
By the It\^o formula, we observe that
 $M^{\epsilon_k}_{N_k,f}(t)$ is a continuous square integrable martingale with quadratic variation process
$$
\big[M^{\epsilon_k}_{N_k,f}\big](t)={1\over N_k}\int_0^t\bigg\langle\mu_{N_k}^{\epsilon_k}(s),\Big(\Big[a\big(\cdot,\mu_{N_k}^{\epsilon_k}(s),\alpha^{\epsilon_k}(s)\big)\nabla_xf(\cdot)\Big]'\nabla_xf(\cdot)\Big)\bigg\rangle ds.
$$
Similar to \eqref{est-quad-var}, since $\sigma(\cdot,\cdot,\cdot)$ is bounded we have
$
\sup_{0\le t\le T}\EE\big[M^{\epsilon_k}_{N_k,f}\big](t)\le {CT\over N_k}.
$
Thus, by Doob's inequality,
\begin{align}
&\PP\bigg(\bigg|\big\langle\mu^{\epsilon_k}_{N_k}(t),f\big\rangle-\big\langle\mu^{\epsilon_k}_{N_k}(0),f\big\rangle-\int_0^t\Big\langle\mu^{\epsilon_k}_{N_k}(s),\mathcal{L}^{\epsilon_k}\big(\mu^{\epsilon_k}_{N_k}(s)\big)f\big(\cdot,\alpha^{\epsilon_k}(s)\big)\Big\rangle ds\bigg|>{\delta\over2}\bigg)\notag\\
&=\PP\bigg(\big|M^{\epsilon_k}_{N_k,f}(t)\big|>{\delta\over2}\bigg) \le {2\over\delta^2}\EE\big[M^{\epsilon_k}_{N_k,f}\big](T) \le {CT\over\delta^2 N_k}.\label{est-prob2-timescale}
\end{align}
Combining \eqref{est-prob-timescale}-\eqref{est-prob2-timescale} yields
$$
\liminf_{k\to\infty}\bar{\mathscr{P}}^{\epsilon_k}_{N_k}\Big(\big|\bar M_f(t)\big|>\delta\Big)=0.
$$
Since $\delta$ is taken arbitrarily, the above equation and \eqref{est-barP1} imply \eqref{suffice-barP}.

Next, we prove that if for some pair $(\eta,\bar\varsigma)$, \eqref{eqn-char-perb} holds  for any $f\in C^2_c\big(\mathbb{R}^d\big)$, it also holds  for any $g(\cdot,\cdot)$ such that $g(\cdot,i)\in C^2_c\big(\mathbb{R}^d\big)$ for each $i\in\bar{\mathbb{S}}$. For $f\in C^2_c\big(\mathbb{R}^d\big)$ and $i\in \bar{\mathbb{S}}$ denote
$$
\hat{\bar{\mathcal{L}}}_i(\mu)f(x)=\bar b'\big(x,\mu,i\big)\nabla_{x}f\big(x\big)+{1\over2}\Big(\bar a(x,\mu,i)\nabla_x\Big)'\nabla_xf\big(x\big).
$$
Then, for any $f\in C^2_c\big(\mathbb{R}^d\big)$ and $g(\cdot,i)\in C^2_c\big(\mathbb{R}^d\big)$, $i\in\bar{\mathbb{S}}$,
\begin{equation}\label{eqn-hatbarL}
\hat{\bar{\mathcal{L}}}_i(\mu)f(x)=\bar{\mathcal{L}}(\mu)f(x,i),\qquad \bar{\mathcal{L}}(\mu)g(x,i)=\hat{\bar{\mathcal{L}}}_i(\mu)g(x,i)+\bar Qg(x)(i).
\end{equation}
Let $(\eta,\bar\varsigma)$ be a pair in $D\big([0,T],\mathscr{M}_1\times\bar{\mathbb{S}}\big)$ such that \eqref{eqn-char-perb} holds  for any $f\in C^2_c\big(\mathbb{R}^d\big)$.  By the definition of $\hat{\bar{\mathcal{L}}}_i$, we can rewrite \eqref{eqn-char-perb} as follow
\begin{equation}\label{eqn-f-short}
\big\langle\eta(t),f\big\rangle=\big\langle\eta(r),f\big\rangle+\int_r^t\Big\langle\eta(s),\hat{\bar{\mathcal{L}}}_{\bar\varsigma(s_-)}\big(\eta(s)\big)f\Big\rangle ds,
\end{equation}
for any $0\le r\le t\le T$ and $f\in C^2_c\big(\mathbb{R}^d\big)$.

Denote the jump times of $\bar\varsigma$ by $t_0=0$, $t_{n+1}=\inf\big\{t>t_n:\bar\varsigma(t)\ne\bar\varsigma(t_-)\big\}$ and $\iota_n=\bar\varsigma(t_n)$ for $n\ge0$. For $t_n<t<t_{n+1}$ we have
\begin{equation}\label{eqn-Mij}
\bar M^{\bar\varsigma}_{ij}(t)-\bar M^{\bar\varsigma}_{ij}(t_n)=
\begin{cases}
&-\Big(\langle \bar M^{\bar\varsigma}_{ij}\rangle(t)-\langle \bar M^{\bar\varsigma}_{ij}\rangle(t_n)\Big)=-q_{\iota_nj}(t-t_n) \quad \text{ if}\quad  j\ne i=\iota_n,\\
&0 \quad \text{ otherwise.}
\end{cases}
\end{equation}
This implies
\begin{equation}\label{eqn-rel-Q-M}
\int_{t_n}^t\Big\langle\eta(s),\bar Qg(\cdot)(\bar\varsigma(s))\Big\rangle ds+\sum_{i,j\in\bar{\mathbb{S}}}\int_{t_n}^t\Big\langle \eta(s),g(\cdot,j)-g(\cdot,i)\Big\rangle d\bar M^{\bar\varsigma}_{ij}(s)=0\quad t_n\le t<t_{n+1}.
\end{equation}
In view of \eqref{eqn-f-short} with $f(x)=g(x,\iota_n)$, \eqref{eqn-hatbarL} and \eqref{eqn-rel-Q-M}, we have
\begin{align}
\big\langle\eta(t),g(\cdot,\bar\varsigma(t))\big\rangle&=\big\langle\eta(t_n),g(\cdot,\bar\varsigma(t_n))\big\rangle+\int_{t_n}^t\Big\langle\eta(s),\hat{\bar{\mathcal{L}}}_{\bar\varsigma(s_-)}\big(\eta(s)\big)g\big(\cdot,\bar\varsigma(s_-)\big)\Big\rangle ds \notag\\
&=\big\langle\eta(t_n),g(\cdot,\bar\varsigma(t_n))\big\rangle+\int_{t_n}^t\Big\langle\eta(s),\bar{\mathcal{L}}\big(\eta(s)\big)g\big(\cdot,\bar\varsigma(s_-)\big)\Big\rangle ds\notag\\
&\quad\quad\quad\quad+\sum_{i,j\in\bar{\mathbb{S}}}\int_{t_n}^t\Big\langle \eta(s),g(\cdot,j)-g(\cdot,i)\Big\rangle d\bar M^{\bar\varsigma}_{ij}(s), \quad t_n\le t<t_{n+1}.\label{eqn-char-perb-1}
\end{align}
At $t=t_{n+1}$, \eqref{eqn-Mij} still holds if $j\ne\iota_{n+1}$.
In addition,
$$\begin{array}{ll}
   &\!\! \disp \big[\bar M^{\bar\varsigma}_{\iota_n\iota_{n+1}}\big](t_{n+1})=
\big[\bar M^{\bar\varsigma}_{\iota_n\iota_{n+1}}\big](t_{n})+1
 \ \hbox{ and }\\
 &\!\! \disp \big[\bar M^{\bar\varsigma}_{ij}\big](t_{n+1})=\big[\bar M^{\bar\varsigma}_{ij}\big](t_{n}) \ \hbox{ if}  \ (i,j)\ne (\iota_n,\iota_{n+1}).
 \end{array}$$
 Thus,
\bea\ad
\int_{t_n}^{t_{n+1}}\Big\langle\eta(s),\bar Qg(\cdot)(\bar\varsigma(s))\Big\rangle ds+\sum_{i,j\in\bar{\mathbb{S}}}\int_{t_n}^{t_{n+1}}\Big\langle \eta(s),g(\cdot,j)-g(\cdot,i)\Big\rangle d\bar M^{\bar\varsigma}_{ij}(s)\\
\aad \ =\Big\langle \eta(t_{n+1}),g(\cdot,\iota_{n+1})-g(\cdot,\iota_n)\Big\rangle.
\eea
As a consequence, by applying  \eqref{eqn-f-short} with $f(x)=g(x,\iota_n)$ and $t=t_{n+1}$ we have
\begin{align*}
&\!\!\!\big\langle\eta(t_{n+1}),g(\cdot,\iota_n)\big\rangle\\
&=\big\langle\eta(t_n),g(\cdot,\iota_n)\big\rangle+\int_{t_n}^{t+1}\Big\langle\eta(s),\hat{\bar{\mathcal{L}}}_{\iota_n}\big(\eta(s)\big)g\big(\cdot,\iota_n\big)\Big\rangle ds+\int_{t_n}^{t_{n+1}}\Big\langle\eta(s),\bar Qg(\cdot)(\iota_n)\Big\rangle ds\\
&\quad\quad\quad\quad +\sum_{i,j\in\bar{\mathbb{S}}}\int_{t_n}^{t_{n+1}}\Big\langle \eta(s),g(\cdot,j)-g(\cdot,i)\Big\rangle d\bar M^{\bar\varsigma}_{ij}(s)-\Big\langle \eta(t_{n+1}),g(\cdot,\iota_{n+1})-g(\cdot,\iota_n)\Big\rangle\\
&=\big\langle\eta(t_n),g(\cdot,\bar\varsigma(t_n))\big\rangle+\int_{t_n}^{t_{n+1}}\Big\langle\eta(s),\bar{\mathcal{L}}\big(\eta(s)\big)g\big(\cdot,\bar\varsigma(s_-)\big)\Big\rangle ds\\
&\quad\quad\quad\quad
+\sum_{i,j\in\bar{\mathbb{S}}}\int_{t_n}^{t_{n+1}}\Big\langle \eta(s),g(\cdot,j)-g(\cdot,i)\Big\rangle d\bar M^{\bar\varsigma}_{ij}(s)\\
&\quad\quad\quad\quad+\Big\langle \eta(t_{n+1}),g(\cdot,\iota_n)-g(\cdot,\iota_{n+1})\Big\rangle.
\end{align*}
This implies that \eqref{eqn-char-perb-1} also holds for $t=t_{n+1}$ and therefore \eqref{eqn-char-perb-1} is proved for any $g(\cdot,\cdot)$ such that $g(\cdot,i)\in C^2_c\big(\mathbb{R}^d\big)$ for each $i\in \bar{\mathbb{S}}$ as desired.
\qed

As a direct consequence of \thmref{thm-SMV} we have the following proposition, which characterizes the limit $\bar{\mathscr{P}}$ as a solution to a stochastic McKean-Vlasov equation with Markovian switching.

\begin{prop}\label{prop-SMV-TwoTimeScale}
Assume {\rm(A1)}, {\rm(A2)}, and {\rm(A3)}. Let $\mu_0$ be a measure in $\mathscr M_1$. Then the system of integral equations
\begin{align}
\big\langle\eta(t),f(\cdot,\bar\alpha(t))\big\rangle=\,\,&\big\langle\mu_0,f(\cdot,\bar\alpha(0))\big\rangle+\int_0^t\Big\langle\eta(s),\bar{\mathcal{L}}\big(\eta(s)\big)f(\cdot,\bar\alpha(s_-))\Big\rangle ds\notag\\
&+\sum_{i,j\in\bar{\mathbb{S}}}\int_0^t\Big\langle\eta(s),f(\cdot,j)-f(\cdot,i)\Big\rangle d{\bar M}_{ij}(s),\,\,0\le t\le T, \label{eq-perturb-SMV-Scale}
\end{align}
where $f(\cdot,i)\in C^2_c(\mathbb{R}^d)$ for each $i\in\bar{\mathbb{S}}$, has a unique solution in $D\big([0,T],\mathscr{M}_1\big)$. Moreover, this solution equals $\mathscr{L}\big(\bar x(t)\big|\mathcal{F}^{\bar\alpha}_{t_-}\big)$ for all $0\le t\le T$,
 where $\bar x(t)$ is the unique solution of
$$
\begin{cases}
d\bar x(t)&=\bar b\Big(\bar x(t),\eta_{\bar\alpha}(t),\bar \alpha(t_-)\Big)dt+\bar \sigma\Big(x(t),\eta_{\bar\alpha}(t),\bar \alpha(t_-)\Big)d\tilde w(t),\quad \mathscr L\big(\bar x(0)\big)=\mu_0,\\
\eta_{\bar \alpha}(t)&=\mathscr{L}\big(\bar x(t)\big|\mathcal{F}^{\bar \alpha}_{t_-}\big)
\end{cases}
$$
where $\tilde w(\cdot)$ is a standard Brownian motion independent of $\bar \alpha(\cdot)$.
\end{prop}

\para{Proof of \thmref{thm-main2}.} Since $\bar\alpha^\epsilon(\cdot)$ converges weakly to $\bar\alpha(\cdot)$, for any $\mathcal{C}\subset D\big([0,T],\bar{\mathbb{S}}\big)$ we have
$\bar{\mathscr{P}}\big(\mathcal{C}\big)=\PP\big(\bar\alpha_T\in \mathcal{C}\big)$. By using $\bar{\mathcal{L}}, \big(\bar\alpha(\cdot),\bar{\mathbb{S}}\big)$, respectively, instead of $\mathcal{L}, \big(\alpha(\cdot), {\mathbb{S}}\big)$, \propref{Prop-SDE-char} instead of \thmref{thm-limit-characterization}, and \propref{prop-SMV-TwoTimeScale} instead of \thmref{thm-SMV}, a similar argument to that in the proof of \thmref{thm-main1} yields the assertion of \thmref{thm-main2}.\qed

\appendix
\section{Appendix}
\subsection{Proof of Lemma \ref{lem-moments}}
 Define
$$
\tau_k=\tau_{k,N}=\inf\Big\{t\ge0:x_i(t)>k \text{ for some }i, 1\le i\le N\Big\}.
$$
The sequence $\tau_k, k=1,2,\ldots$ is monotonically increasing. Put $\tau_\infty=\lim_{k\to\infty}\tau_k$.

We are in a position to prove that $\lim_{k\to\infty}\tau_k=\infty$ a.s. Suppose that there exists a positive numbers $T_0$ and $\epsilon$ such that $P\big(\lim_{k\to\infty}\tau_k<T_0\big)>2\epsilon$. Then there is a number $k_0$ such that $P\big(\tau_k<T_0\big)>\epsilon$ for all $k\ge k_0$. Recall that  for $x\in\mathbb{R}^d$, $\psi(x)=|x|^2$. For $\underline{x}=(x_1,x_2,\ldots, x_N)\in (\mathbb{R}^d)^N$, define the Lyapunov function
$$
V(\underline{x},\alpha)=\bigg({1\over N}\sum_{i=1}^N|x_i|^2+1\bigg)^p.
$$
It is easily seen that
\begin{equation}\label{nablaV}
\nabla_{x_i}V(\underline{x},\alpha)={2p\over N}\bigg(\sum_{i=1}^N|x_i|^2+1\bigg)^{p-1}x_i',
\end{equation}
\begin{equation}\label{HessianV}
\nabla_{x_i}^2V(\underline{x},\alpha)={2p\over N}\bigg(\sum_{i=1}^N|x_i|^2+1\bigg)^{p-1}I_{d\times d}+{4p(p-1)\over N^2}\bigg(\sum_{i=1}^N|x_i|^2+1\bigg)^{p-2}x_ix_i',
\end{equation}
where $\nabla_{x_i}^2V=[\nabla_{x_i}(\nabla_{x_i}V)]'$ denotes the $d\times d$ Hessian matrix with respect to the variable $x_i$ of V.  Since $w_1(\cdot),w_2(\cdot),\ldots, w_N(\cdot)$ are independent Brownian motions, a direct calculation yields
\begin{equation}\label{L_N(V)}
\begin{array}{ll}\disp
\mathcal L_NV(\underline{x},\alpha) &\!\!\! \disp =\sum_{i=1}^N\nabla_{x_i}V(\underline{x},\alpha)b\Big(x_i,{1\over N}\sum_{j=1}^N\delta_{x_j},\alpha\Big)\\
&\quad \  \disp +{1\over2}
\sum_{i=1}^N\text{trace}\Bigg[\nabla_{x_i}^2V\big(\underline{x},\alpha\big)a\Big(x_i,{1\over N}\sum_{j=1}^N\delta_{x_j},\alpha\Big)\Bigg].\end{array}
\end{equation}
It follows from Assumption A and equations \eqref{nablaV}-\eqref{L_N(V)} that
\begin{equation}\label{L_N(V)-ineq}
\big|\mathcal L_NV\big(\underline{x},\alpha\big)\big|\le CV(\underline{x},\alpha),
\end{equation}
for some constant $C$ independent of $N$. Denote $\underline{x}(t)=\big(x_1(t),x_2(t),\ldots,x_N(t)\big)'$. Then the Dynkin formula implies
\begin{align*}
\EE V\big(\underline{x}(T_0\wedge\tau_k),\alpha(T_0\wedge\tau_k)\big)&=\EE V\big(\underline{x}(0),\alpha(0)\big)+\int_0^{T_0\wedge\tau_k}\EE\mathcal L_N V(\underline{x}(s),\alpha(s))ds\\
&\le \EE V\big(\underline{x}(0),\alpha(0)\big)+C\int_0^{T_0\wedge\tau_k}
\EE  V(\underline{x}(s),\alpha(s))ds.
\end{align*}
By the Gronwall inequality, we obtain
\begin{equation}\label{EV-inq}
\EE V\big(\underline{x}(T_0\wedge\tau_k),\alpha(T_0\wedge\tau_k)\big)\le e^{CT_0}\EE V\big(\underline{x}(0),\alpha(0)\big).
\end{equation}
According to the definitions of $\tau_k$ and $V$,
$$
V\big(\underline{x}(T_0\wedge\tau_k),\alpha(T_0\wedge\tau_k)\big)I_{\{\tau_k<T_0\}}=V\big(\underline{x}(\tau_k),\alpha(\tau_k)\big)I_{\{\tau_k<T_0\}}\ge \big(k^2+1\big)^pI_{\{\tau_k<T_0\}}.
$$
Thus, \eqref{EV-inq} yields
\begin{align*}
\big(k^2+1\big)^p\PP\Big(\tau_k<T_0\Big)&\le \EE\Big[V\big(\underline{x}(T_0\wedge\tau_k),\alpha(T_0\wedge\tau_k)\big)I_{\{\tau_k<T_0\}}\Big]\\
&\le \EE V\big(\underline{x}(T_0\wedge\tau_k),\alpha(T_0\wedge\tau_k)\big)\\
&\le e^{CT_0}\EE V\big(\underline{x}(0),\alpha(0)\big).
\end{align*}
As a consequence,
$$
\epsilon< \PP\Big(\tau_k<T_0\Big)\le \big(k^2+1\big)^{-p}e^{CT_0}\EE V\big(\underline{x}(0),\alpha(0)\big)\quad \forall k\ge k_0.
$$
This is a contradiction. As a result, $\lim_{k\to\infty}\tau_k=\infty$ a.s.

Next, by applying \eqref{EV-inq} for $t$ instead of $T_0$, with $0\le t\le T$, and letting $k\to\infty$, we arrive at
$$
\EE V\big(\underline{x}(t),\alpha(t)\big)\le e^{Ct}\EE V\big(\underline{x}(0),\alpha(0)\big)\le e^{CT}\EE V\big(\underline{x}(0),\alpha(0)\big).
$$
Since $V\big(\underline{x}(t),\alpha(t)\big)=\big[\langle \mu_N(t),\psi\rangle+1\big]^p$, this inequality implies \eqref{moment-inq} as desired.

To proceed, using the Dynkin formula again, we obtain
\begin{align*}
&\EE\Big(V\big(\underline{x}(t\wedge\tau_k),\alpha(t\wedge\tau_k)\big)\Big|{\mathcal F}^{N,\alpha}_{s\wedge\tau_k}\Big)\\
&=V\big(\underline{x}(s\wedge\tau_k),\alpha(s\wedge\tau_k)\big)+\int_{s\wedge\tau_k}^{t\wedge\tau_k}\EE\Big(\mathcal L_NV\big(\underline{x}(u),\alpha(u)\big)\Big|{\mathcal F}^{N,\alpha}_{s\wedge\tau_k}\Big)du
\end{align*}
for $0\le s\le t$. It follows from \eqref{L_N(V)-ineq} that
\begin{align*}
&V\big(\underline{x}(s\wedge\tau_k),\alpha(s\wedge\tau_k)\big)-C\int_{s\wedge\tau_k}^{t\wedge\tau_k}\EE\Big(V\big(\underline{x}(u\wedge\tau_k),\alpha(u\wedge\tau_k)\big)\Big|{\mathcal F}^{N,\alpha}_{s\wedge\tau_k}\Big)du\\
&\quad \le \EE\Big(V\big(\underline{x}(t\wedge\tau_k),\alpha(t\wedge\tau_k)\big)\Big|{\mathcal F}^{N,\alpha}_{s\wedge\tau_k}\Big)\\
&\quad \le V\big(\underline{x}(s\wedge\tau_k),\alpha(s\wedge\tau_k)\big)+C\int_{s\wedge\tau_k}^{t\wedge\tau_k}\EE\Big(V\big(\underline{x}(u\wedge\tau_k),\alpha(u\wedge\tau_k)\big)\Big|{\mathcal F}^{N,\alpha}_{s\wedge\tau_k}\Big)du.
\end{align*}
Letting $k\to\infty$, by virtue of the Gronwall inequality, there exists a constant $C$ independent of $N$ such that for $0\le s\le t$,
$$
e^{-C(t-s)}\EE \Big(V\big(\underline{x}(s),\alpha(s)\big)\Big|{\mathcal F}^{N,\alpha}_s\Big)\le \EE\Big(V\big(\underline{x}(t),\alpha(t)\big)\Big|{\mathcal F}^{N,\alpha}_s\Big)\le e^{C(t-s)}\EE\Big( V\big(\underline{x}(s),\alpha(s)\big)\Big|{\mathcal F}^{N,\alpha}_s\Big).
$$
This gives \eqref{cond-moment-inq} and completes the proof.
\qed

\subsection{Proof of Lemma \ref{lem-Lambda}}
We prove the two parts of the assertions as follows.

(i) Let $\varsigma_1,\varsigma_2\in D_f\big([0,r),\mathbb{S}\big)$ and denote $\eta_1=\Lambda_r(\varsigma_1), \eta_2=\Lambda_r(\varsigma_2)$. It follows from Step 1 in the proof of \thmref{thm-SMV} that for $i=1,2$, $\eta_i=\mathscr{L}(y_i(s))$  which are the unique solution to the equation
\begin{equation}\label{eqn-yi}
dy_i(s)=b\Big(y_i(s),\eta_i(s),\varsigma_i(s_-)\Big)ds+\sigma\Big(y_i(s),\eta_i(s),\varsigma_i(s_-)\Big)d\tilde w(s),\quad 0\le s\le r,
\end{equation}
where $\mathscr L\big(y_i(0)\big)=\mu_0$. Without loss of generality, we can assume that $y_1(0)=y_2(0)$.

First, we show that there is a constant $C$ independent on $\varsigma_i$ and $r$ such that
\begin{equation}\label{inq-moment-yi}
\EE\Big(\sup_{0\le s\le r}\big|y_i(s)\big|^2\Big)\le C.
\end{equation}
It follows from \eqref{eqn-yi} that
\begin{equation}\label{inq-|yi|^2}
|y_i(t)|^2\le 3|y_i(0)|^2+3T\int_0^t\Big|b\Big(y_i(s),\eta_i(s),\varsigma_i(s_-)\Big)\Big|^2ds+3\Big|\int_0^t\sigma\Big(y_i(s),\eta_i(s),\varsigma_i(s_-)\Big)d\tilde w(s)\Big|^2.
\end{equation}
Note that $\langle\eta_i(s),\varphi\rangle^2=\big(E|y_i(s)|\big)^2\le E|y_i(s)|^2$. Thus, by  Assumption A, there is a constant $C$ independent on $\varsigma_i$ and $r$ such that for each $0\le s\le r$, $\sigma\big(y_i(s),\eta_i(s),\varsigma_i(s_-)\big)\le C$ and
\begin{equation}\label{inq-b(yi)}
\big|b\big(y_i(s),\eta_i(s),\varsigma_i(s_-)\big)\big|^2\le C\big(1+|y_i(s)|^2+\EE|y_i(s)|^2\big).
\end{equation}
Taking these inequalities into account and using the Burkholder-Davis-Gundy inequality for the last term in the right-hand side of \eqref{inq-|yi|^2}, we arrive at
$$
\sup_{0\le s\le t}|y_i(s)|^2\le C+C\int_0^t\Big[\sup_{0\le u\le s}|y_i(u)|^2+\EE\big(\sup_{0\le u\le s}|y_i(u)|^2\big)\Big]ds.
$$
By taking expectations on both sides of the above inequality and using the Gronwall inequality we obtain \eqref{inq-moment-yi}.

In view of \eqref{inq-moment-yi} and \eqref{inq-b(yi)}, by assumption {\rm (A2)} we have
$$
\big\|\eta_i(t+\delta)-\eta_i(t)\big\|_{BL}^2\le E\big|y_i(t+\delta)-y_i(t)\big|^2\le C\delta,\quad i=1,2,
$$
for any $\delta>0$ such that $0\le t<t+\delta<r$ where $C$ is independent of $\delta$. This shows that $\eta_1,\eta_2\in C\big([0,r),\mathscr{M}_1\big)$.

Next, we prove the continuity of $\Lambda_r$. Let $\rho$ be the metric on $D([0,r),\mathbb{S})$ defined by
$$
\rho(\varsigma_1,\varsigma_2)=\inf_{\lambda}\Big\{\sup_{t\in[0,r)}d_{\mathbb{S}}\big(\varsigma_1(t),\varsigma_2(\lambda(t))\big),\sup_{t\in[0,r)}\big|t-\lambda(t)\big|\Big\},\quad \varsigma_1,\varsigma_2\in D([0,r),\mathbb{S}),
$$
where the infimum is taken over the class of all strictly increasing and continuous mappings $\lambda$ of $[0,r)$ onto itself. Because of the definition of the metric $d_{\mathbb{S}}$ on the discrete set $\mathbb{S}$, if $\rho(\varsigma_1,\varsigma_2)<1$, $\varsigma_1$ and $\varsigma_2$ have same number of jumps on $[0,r)$.
Hence, to prove the continuity of $\Lambda_r$, it suffices to show that there is a constant $C$ such that
\begin{equation}\label{inq-suff-cts}
\sup_{0\le t<r}\|\eta_1(t)-\eta_2(t)\|_{BL}^2\le Cm\rho(\varsigma_1,\varsigma_2),\quad \varsigma_1,\varsigma_2\in D_f([0,r),\mathbb{S}),
\end{equation}
whenever $\delta=\rho(\varsigma_1,\varsigma_2)<1$ where $m$ is the number of jumps on $[0,r)$ of $\varsigma_1$ and $\varsigma_2$. According to the definition of $\|\cdot\|_{BL}$ and $\big\langle\eta_1(t)-\eta_2(t),f\rangle=\EE\big(f(y_1(t))-f(y_2(t))\big)$, we have
\begin{align}
\big\|\eta_1(t)-\eta_2(t)\big\|_{BL}\le \EE\big|y_1(t)-y_2(t)\big|.\label{inq-eta12}
\end{align}
Since $d_{\mathbb{S}}(i_0,j_0)=1$ if $i_0\ne j_0$, we have
$$
I_1=\{t\in[0,r):\varsigma_1(t)\ne\varsigma_2(t)\}=\cup_{k=1}^mA_k,
$$
where $m$ is the number of jumps of $\varsigma_1$ and $\varsigma_2$, and each $A_k$ is an interval with length $|A_k|<\delta$. Denote $I_0=[0,r)\backslash I_1$ and
\begin{align}
J_i(t)=&\EE\int_{I_i\cap[0,t]}\Big|b\big(y_1(s),\eta_1(s),\varsigma_1(s_-)\big)-b\big(y_2(s),\eta_2(s),\varsigma_2(s_-)\big)\Big|^2ds\notag\\
&\quad+\EE\Big|\int_{I_i\cap[0,t]}\Big[\sigma\big(y_1(s),\eta_1(s),\varsigma_1(s_-)\big)-\sigma\big(y_2(s),\eta_2(s),\varsigma_2(s_-)\big)\Big]d\tilde w(s)\Big|^2,\quad i=0,1.\label{def-J01}
\end{align}
According to \eqref{eqn-yi} and the Cauchy-Schwarz inequality, there is a constant $C=C(T)$ such that
\begin{equation}\label{inq-J01}
\EE\big|y_1(t)-y_2(t)\big|^2\le C\big(J_0(t)+J_1(t)\big).
\end{equation}
To proceed, we estimate $J_0(t)$ and $J_1(t)$. Since $\varsigma_1=\varsigma_2$ on $I_0$, assumption (A1) and \eqref{inq-eta12} give
\begin{align}
&\Big|b\big(y_1(s),\eta_1(s),\varsigma_1(s_-)\big)-b\big(y_2(s),\eta_2(s),\varsigma_2(s_-)\big)\Big|^2\notag\\
&\hspace{5cm}+\Big|\sigma\big(y_1(s),\eta_1(s),\varsigma_1(s_-)\big)-\sigma\big(y_2(s),\eta_2(s),\varsigma_2(s_-)\big)\Big|^2\notag\\
&\le2\Big(\big|y_1(s)-y_2(s)\big|^2+\EE\big|y_1(s)-y_2(s)\big|^2\Big), \label{inq-|b1-b2|}
\end{align}
for $s\in I_0$ excepts at finite points on its boundary. Therefore, by virtue of the Burkholder-Davis-Gundy inequality, \eqref{def-J01} and \eqref{inq-|b1-b2|} imply
\begin{equation}\label{inq-J0}
J_0(t)\le C\int_0^t\EE\big|y_1(s)-y_2(s)\big|^2ds,
\end{equation}
where $C=C(T)$ is a constant which only depends on $T$.
Now we are in a position to estimate $J_1(t)$. It follows from \eqref{inq-moment-yi} and \eqref{inq-b(yi)} that
$$
\EE\Big(\sup_{0\le s\le r}\Big|b\big(y_1(s),\eta_1(s),\varsigma_1(s_-)\big)-b\big(y_2(s),\eta_2(s),\varsigma_2(s_-)\big)\Big|^2\Big)\le C.
$$
Since $I_1=\cup_{k=1}^mA_k$ where $m<\infty$, by the boundedness of $\sigma(\cdot,\cdot,\cdot)$, and the Cauchy-Schwarz and the Burkholder-Davis-Gundy  inequalities,
\begin{align}
J_1(t)&\le C\sum_{i=1}^m\int_{A_i}\EE\Big(\sup_{0\le s\le r}\Big|b\big(y_1(s),\eta_1(s),\varsigma_1(s_-)\big)-b\big(y_2(s),\eta_2(s),\varsigma_2(s_-)\big)\Big|^2\Big)ds\notag\\
&\quad+C\sum_{i=1}^m\Big|\int_{A_i}\Big[\sigma\big(y_1(s),\eta_1(s),\varsigma_1(s_-)\big)-\sigma\big(y_2(s),\eta_2(s),\varsigma_2(s_-)\big)\Big]d\tilde w(s)\Big|^2\notag\\
&\le C\sum_{i=1}^m|A_i|\le Cm\delta.\label{inq-J1}
\end{align}
Combining \eqref{inq-J01}, \eqref{inq-J0}, and \eqref{inq-J1}, we obtain
$$
\EE\big|y_1(t)-y_2(t)\big|^2\le C\int_0^tE\big|y_1(s)-y_2(s)\big|^2ds+Cm\delta.
$$
By the Gronwall inequality, for any $0\le t<r$,
$$\EE\big|y_1(t)-y_2(t)\big|^2\le Cm\delta e^{Cr},$$ which together with \eqref{inq-eta12} implies \eqref{inq-suff-cts}.

 (ii) Equation
 \eqref{eqn-consistency} follows from the uniqueness of the solution of \eqref{eq-SMV-varsig} proved in Step 1 in the proof of \thmref{thm-SMV}.

\subsection{Proof of Lemma \ref{Lem-operator-est}} We have
\begin{align}
&\int_0^t\Big\langle\mu_N^\epsilon(s),\mathcal{L}^\epsilon\big(\mu_N^\epsilon(s)\big)f(\alpha^\epsilon(s))\Big\rangle ds-\int_0^t\Big\langle\mu_N^\epsilon(s),\bar{\mathcal{L}}(\mu_N^\epsilon(s))f(\bar\alpha^\epsilon(s))\Big\rangle ds\notag\\
=&\int_0^t\bigg\langle\mu_N^\epsilon(s),\sum_{i=1}^{l_0}\sum_{j=1}^{m_i}b'(x,\mu_N^\epsilon(s),s_{ij})\nabla_xf(x)\Big[I(\alpha^\epsilon(s)=s_{ij})-\nu_{ij}I(\bar\alpha^\epsilon(s)=i)\Big]\bigg\rangle ds\notag\\
&+{1\over2}\int_0^t\bigg\langle\mu_N^\epsilon(s),\sum_{i=1}^{l_0}\sum_{j=1}^{m_i}\big(a(x,\mu_N^\epsilon(s),s_{ij})\nabla_x\big)'\nabla_xf(x)\Big[I\big(\alpha^\epsilon(s)=s_{ij}\big)-\nu_{ij}I\big(\bar\alpha^\epsilon(s)=i\big)\Big]\bigg\rangle ds\notag\\
=&{1\over N}\sum_{k=1}^N\sum_{i=1}^{l_0}\sum_{j=1}^{m_i} \int_0^tb'\big(x_k^\epsilon(s),\mu_N^\epsilon(s),s_{ij}\big)\nabla_xf\big(x_k^\epsilon(s)\big)\Big[I\big(\alpha^\epsilon(s)=s_{ij}\big)-\nu_{ij}I\big(\bar\alpha^\epsilon(s)=i\big)\Big] ds\notag\\
&+{1\over 2 N}\sum_{k=1}^N\sum_{i=1}^{l_0}\sum_{j=1}^{m_i}
\int_0^t\Big(a\big(x_k^\epsilon(s),\mu_N^\epsilon(s),s_{ij}\big)
\nabla_x\Big)'\nabla_xf\big(x_k^\epsilon(s)\big)\notag\\
&\qquad \times \Big[I\big(\alpha^\epsilon(s)=s_{ij}\big)-\nu_{ij}I\big(\bar\alpha^\epsilon(s)=i\big)\Big] ds\notag\\
=:&{1\over N}\sum_{k=1}^N\sum_{i=1}^{l_0}\sum_{j=1}^{m_i}I_1(t,k,i,j)+{1\over 2 N}\sum_{k=1}^N\sum_{i=1}^{l_0}\sum_{j=1}^{m_i}I_2(t,k,i,j).
\end{align}
Thus, to prove \eqref{operator-est}, it suffices to show that
\begin{equation}
\sup_{0\le t\le T}\EE\Big[\big|I_1(t,k,i,j)\big|+\big|I_2(t,k,i,j)\big|\Big]\le C\epsilon^{1/6}, \label{I1+I2}
\end{equation}
for some constant $C$ independent of $N$, $i$, $j$, and $k$. Denote
$$
n=\lfloor \epsilon^{-1/3}\rfloor,\
h=T/n=O(\epsilon^{1/3}), \ \hbox{ and } \ t_l=t_{l,n}=lh \ \hbox{ for }\  l=0,1,\ldots,n,
$$
where $\lfloor \epsilon^{-1/3}\rfloor$ is the  greatest integer that is less than or equal to $\epsilon^{-1/3}$. Then $t_l-t_{l-1}=O(\epsilon^{1/3})$. To proceed, we estimate $\big|I_1(t,k,i,j)\big|$. Let $q$ be an integer $0\le q\le n-1$ such that $t_q\le t< t_{q+1}$. Note that
\begin{align}
&\!\!\!\!\!\!\big|I_1(t_q,k,i,j)\big|\notag\\
\le&\sum_{l=1}^n\bigg|\int_{t_{l-1}}^{t_l}b'\big(x_k^\epsilon(s),\mu_N^\epsilon(s),s_{ij}\big)\nabla_xf\big(x^\epsilon_k(s)\big)\Big[I\big(\alpha^\epsilon(s)=s_{ij}\big)-\nu_{ij}I\big(\bar\alpha^\epsilon(s)=i\big)\Big]ds\bigg|\notag\\
\le&2\sum_{l=1}^n\int_{t_{l-1}}^{t_l}\Big|b'\big(x_k^\epsilon(s),\mu_N^\epsilon(s),s_{ij}\big)\nabla_xf\big(x^\epsilon_k(s)\big)-b'\big(x_k^\epsilon(t_{l-1}),\mu_N^\epsilon(t_{l-1}),s_{ij}\big)\nabla_xf(x^\epsilon_k(t_{l-1}))\Big|ds\notag\\
&+\sum_{l=1}^n\Big|b'\big(x_k^\epsilon(t_{l-1}),\mu_N^\epsilon(t_{l-1}),s_{ij}\big)\nabla_xf\big(x^\epsilon_k(t_{l-1})\big)\Big|\bigg|\int_{t_{l-1}}^{t_l}\Big[I\big(\alpha^\epsilon(s)=s_{ij}\big)-\nu_{ij}I\big(\bar\alpha^\epsilon(s)=i\big)\Big]ds\bigg|\notag\\
=:&\sum_{l=1}^n\Big(2J_1(l)+J_2(l)\Big).\label{J1+J2}
\end{align}
We have
\begin{align}
\big|J_1(l)\big|\le & \int_{t_{l-1}}^{t_l}\bigg|\Big[b'\big(x_k^\epsilon(s),\mu_N^\epsilon(s),s_{ij}\big)-b'\big(x_k^\epsilon(t_{l-1}),\mu_N^\epsilon(t_{l-1}),s_{ij}\big)\Big]\nabla_xf\big(x^\epsilon_k(t_{l-1})\big)\bigg|ds\notag\\
&+\int_{t_{l-1}}^{t_l}\bigg|b'\big(x_k^\epsilon(t_{l-1}),\mu_N^\epsilon(t_{l-1}),s_{ij}\big)\Big[\nabla_xf\big(x^\epsilon_k(s)\big)-\nabla_xf\big(x^\epsilon_k(t_{l-1})\big)\Big]\bigg|ds.\notag
\end{align}
Hence, the Cauchy-Schwarz inequality implies
\begin{align}
&\!\!\!\!\!\EE\big|J_1(l)\big|\notag\\
\le & \int_{t_{l-1}}^{t_l}\bigg[\EE\Big|b'\big(x_k^\epsilon(s),
\mu_N^\epsilon(s),s_{ij}\big)-b'\big(x_k^\epsilon(t_{l-1}),
\mu_N^\epsilon(t_{l-1}),s_{ij}\big)\Big|^2\bigg]^{1/2}
\bigg[E\Big|\nabla_xf(x^\epsilon_k(t_{l-1}))\Big|^2\bigg]^{1/2}ds\notag\\
&+\int_{t_{l-1}}^{t_l}\bigg[\EE\Big|b'\big(x_k^\epsilon(t_{l-1}),\mu_N^\epsilon(t_{l-1}),s_{ij}\big)\Big|^2\bigg]^{1/2}\bigg[\EE\Big|\nabla_xf(x^\epsilon_k(s))-\nabla_xf(x^\epsilon_k(t_{l-1}))\Big|^2\bigg]^{1/2}ds.\notag
\end{align}
Since $\nabla_xf(\cdot)$ is Lipschitz and
$$
\sup_{0\le t\le T}\EE\big|b'\big(x_k^\epsilon(t),\mu_N^\epsilon(t),s_{ij}\big)\big|^2<C,
$$
for some constant $C$ independent of $N,i,j$,  assumption  A and \eqref{inq-mu-diff} imply that
\begin{equation}
\EE\big|J_1(l)\big|\le C\int_{t_{l-1}}^{t_l}\epsilon^{1/6}ds=C\epsilon^{1/6}(t_l-t_{l-1}).\label{est-J_1}
\end{equation}
Next, using the Cauchy-Schwarz inequality again, we obtain
\begin{align}
\Big[\EE\big|J_2(l)\big|\Big]^2 & \ \le \EE\Big|b'\big(x_k^\epsilon(t_{l-1}),\mu_N^\epsilon(t_{l-1}),
s_{ij}\big)\nabla_xf(x^\epsilon_k(t_{l-1}))\Big|^2\\
& \qquad\qquad\qquad\qquad \times \EE\bigg|\int_{t_{l-1}}^{t_l}\Big[I\big(\alpha^\epsilon(s)=s_{ij}\big)-\nu_{ij}I\big(\bar\alpha^\epsilon(s)=i\big)\Big]ds\bigg|^2\notag\\
& \ =C\epsilon,\notag
\end{align}
where the last equation is a consequence of Theorem 5.25 of \cite{YinZhang13}.
Therefore,
\begin{equation}
\EE\big|J_2(l)\big|\le C\epsilon^{1/2}.\label{est-J_2}
\end{equation}
Combining \eqref{J1+J2}, \eqref{est-J_1}, and \eqref{est-J_2}, we arrive at
\begin{equation}
\max_{0\le q\le n}\EE\Big|I_1(t_q,k,i,j)\Big|\le C\Big(\epsilon^{1/6}+\epsilon^{1/2-1/3}\Big)=C\epsilon^{1/6}.\label{est-I1}
\end{equation}
By the same argument, we can show that
\begin{align}
&\max_{0\le q\le n}\sup_{t_q\le t<t_{q+1}}\EE\Big|I_1(t,k,i,j)-I_1(t_q,k,i,j)\Big|\notag\\
&\le \max_{0\le q\le n}\sup_{t_q\le t<t_{q+1}}\EE \bigg|\int_{t_{q}}^{t}b'\big(x_k^\epsilon(s),\mu_N^\epsilon(s),s_{ij}\big)\nabla_xf\big(x^\epsilon_k(s)\big)\notag\\
&\qquad\qquad\qquad\qquad\qquad\qquad\times\Big[I\big(\alpha^\epsilon(s)=s_{ij}\big)-\nu_{ij}I\big(\bar\alpha^\epsilon(s)=i\big)\Big]ds\bigg|\notag\\
&\le C\epsilon^{1/2}.\label{inq-last}
\end{align}
Combining \eqref{est-I1} and \eqref{inq-last} yields $\sup_{0\le t\le T}\EE\Big|I_1(t,k,i,j)\Big|\le C\epsilon^{1/6}$. Likewise, because $f\in C^3_c(\mathbb{R}^d)$, we obtain $\EE\Big|I_2(t,k,i,j)\Big|\le C\epsilon^{1/6}$. Thus \eqref{I1+I2} holds. The proof of the lemma is thus complete.
\qed


\begin{thebibliography}{99}
{\small
\setlength{\baselineskip}{0.12in}
\parskip=0pt

\bibitem{LDF17}
 Andreis L.,  Dai Pra P., and  Fischer M.  McKean-Vlasov limit for interacting systems with simultaneous jumps.  {\sl arXiv preprint} arXiv:1704.01052 (2017).

\bibitem{BaladronFFT12}
 Baladron,  J., Fasoli, D., Faugeras, O., and  Touboul J.,  Mean-field description and propagation of chaos in networks of Hodgkin-Huxley and FitzHugh-Nagumo neurons. {\sl J. Mathematical Neuroscience} {\bf2} (2012), p. 10.



 \bibitem{BSY16}
 Bao, J., Shao, J. and Yuan, C. Approximation of invariant measures for regime-switching diffusions, {\sl Potential Anal.} {\bf 44} (2016), 707--727.

\bibitem{BFY13}
Bensoussan, A., Frehse, J. and Yam, P. {\it Mean Field Games and Mean Field Type Control Theory}, SpringerBriefs in Mathematics. Springer, New York, 2013.

\bibitem{ContucciGM08}
Contucci, P., Gallo,  I. and  Menconi, G. Phase transitions in social sciences: two-populations mean field theory.  {\sl International Journal of Modern Physics B} {\bf 22} (2008), 2199--2212.

\bibitem{CD13} Costa, O. and Dufour, F.
{\it Continuous Average Control of Piecewise Deterministic Markov Processes}, Springer, New York, 2013.




\bibitem{Dawson83}
Dawson, D.A. Critical dynamics and fluctuations for a mean-field model of cooperative behavior. {\sl J. Statist. Phys.} {\bf31} (1983), no. 1, 29--85.





\bibitem{DawsonV95}
Dawson, D. and Vaillancourt, J. Stochastic McKean-Vlasov equations. {\sl NoDEA Nonlinear Differential Equations Appl.} {\bf2} (1995), no. 2, 199--229.

\bibitem{DZ1991} Dawson, D.A.  and Zheng, X. Law of large numbers and central
limit theorem for unbounded jump mean-field models, {\it Adv. Appl.
Math.}, 12 (1991), 293--326.

\bibitem{Dudley66}
Dudley, R.M. Convergence of Baire measures. {\sl Studia Math.} {\bf 27} (1966), 251--268.

\bibitem{EthierKurtz86}
Ethier, S.N. and Kurtz, T.G. {\sl Markov Processes. Characterization and Convergence.} Wiley Series in Probability and Mathematical Statistics. John Wiley \& Sons, Inc., New York, 1986.

\bibitem{FlemingR}
Fleming, W.H. and  Rishel, R.W.
{\it Deterministic and Stochastic Optimal Control}, Springer-Verlag,
New York, NY, 1975.


\bibitem{Gartner88}
G\"artner, J. On the McKean-Vlasov limit for interacting diffusions. {\sl Math. Nachr.} {\bf137} (1988), 197--248.

\bibitem{Graham90}
Graham, C.  Nonlinear limit for a system of diffusing particles which alternate between two states. {\sl Appl. Math. Optim.} {\bf 22} (1990), 75--90.


\bibitem{HN16}
 Huang, M. and Nguyen, S. Mean field games for stochastic growth with relative utility, {\sl Appl. Math. Optim.} {\bf 74} (2016),  643--668.

\bibitem{HCM03}
 Huang, M.,  Caines P.E. and  Malham\'e, R.P. Individual and mass behavior in large population stochastic wireless power
control problems: centralized and Nash equilibrium solutions. (2003) 98--103.


\bibitem{HMC06}
 Huang, M.,  Malham\'e, R.P., and  Caines, P.E. Large population stochastic
dynamic games: Closed-loop McKean-Vlasov systems and the Nash
certainty equivalence principle. {\sl Comm. Information Syst.} {\bf6} (2006), 221--252.

\bibitem{Kac56}
Kac, M. Foundations of kinetic theory. In {\sl Proceedings of the Third Berkeley Symposium on Mathematical Statistics and Probability}, Univ of California Press, {\bf3} (1956), 171--197.

\bibitem{Kolokoltsov10}
Kolokoltsov, V. N. {\sl Nonlinear Markov Processes and Kinetic Equations}, volume 182. Cambridge University
Press, 2010.


\bibitem{KurtzXiong99}
Kurtz, T.G. and Xiong, J. Particle representations for a class of nonlinear SPDEs. {\sl Stochastic Process. Appl.} {\bf83} (1999), no. 1, 103--126.


\bibitem{LL06a}
 Lasry, J.M. and  Lions, P.L. Jeux \`a champ moyen. I. Le cas stationnaire. {\sl C.
R. Math. Acad. Sci. Paris} {\bf343} (2006) 619--625.



\bibitem{LiptserShiryayev89}
Liptser, R.Sh. and Shiryayev, A.N. {\sl Theory of Martingales}.   Kluwer Academic Publishers Group, Dordrecht, 1989.


\bibitem{LM09}
Luo, Q. and Mao, X. Stochastic population dynamics under regime switching, II, {\sl J. Math. Anal. Appl.} {\bf 355} (2009),  577--593.

\bibitem{MaoYuan06}
Mao, X. and Yuan, C. {\sl Stochastic Differential Equations with Markovian Switching}. Imperial College Press, London, 2006.



\bibitem{McKean66}
 McKean, H.P. A class of Markov processes associated with nonlinear
parabolic equations. {\sl Proc. Nat. Acad. Sci. U.S.A.}, {\bf56} (1966), 1907--1911.


\bibitem{NH12}
 Nguyen, S. and Huang, M. Linear-quadratic-Gaussian mixed games with continuum-parametrized minor players, {\sl SIAM J. Control Optim.} {\bf 50} (2012), 2907--2937.

\bibitem{Oelschlager84}
Oelschlager, K. A martingale approach to the law of large numbers for weakly interacting stochastic processes. {\sl Ann. Probab.} {\bf12} (1984), no. 2, 458--479.

\bibitem{RogersWilliams00}
Rogers, L.C.G. and  Williams  D. {\sl Diffusions, Markov Processes, and Martingales. Vol. 2. It\^o Calculus}. Cambridge University Press, Cambridge, 2000.

\bibitem{SethiZ}
Sethi, S.P. and   Zhang, Q.
{\sl Hierarchical Decision Making in Stochastic Manufacturing
Systems}, Birkh\"auser, Boston, 1994.



\bibitem{Sznitman91}
 Sznitman, A.S. Topics in propagation of chaos. In P.-L. Hennequin, editor,
Ecole d'Et\'e de Probabilit\'es de Saint-Flour XIX - 1989,
{\sl Lecture Notes in Math.}, {\bf 1464} (1991), 165--251, Springer-Verlag, Berlin.


\bibitem{TYW16}
Tran, K.,  Yin, G. and Wang, L.Y.
A generalized Goodwin business cycle model in random environment,
 {\sl J. Math. Anal. Appl.}, {\bf 438} (2016), 311--327.

\bibitem{WangZ12}
Wang, B.C. and  Zhang, J.F.  Mean field games for large-population multiagent systems with Markov jump parameters.  {\sl SIAM Journal on Control and Optimization} {\bf50} (2012), 2308--2334.

\bibitem{WangZ17}
Wang, B.C. and  Zhang, J.F. Social Optima in Mean Field Linear-Quadratic-Gaussian Models with Markov Jump Parameters. {\sl SIAM Journal on Control and Optimization} {\bf 55} (2017), 429--456.

\bibitem{XiYin09}
Xi, F. and Yin, G. Asymptotic properties of a mean-field model with a continuous-state-dependent switching process. {\sl J. Appl. Probab.} {\bf46} (2009), no. 1, 221--243.

\bibitem{YKI04}
 Yin, G.,  Krishnamurthy, V. and Ion, C.
Regime switching stochastic approximation algorithms with
application to adaptive discrete stochastic optimization,
{\sl SIAM J. Optim.}, {\bf 14} (2004), 1187--1215.




\bibitem{YinZhang13}
Yin, G. and Zhang, Q. {\sl Continuous-Time Markov Chains and Applications: A Two-Time-Scale Approach}, 2nd Ed.,  Springer, New York, 2013.

\bibitem{YinZhu09}
Yin, G.  and Zhu, C. {\sl Hybrid Switching Diffusions. Properties and Applications},
Springer, New York, 2010.

\bibitem{ZhY03}
 Zhou, X.Y. and  Yin, G.
Markowitz mean-variance portfolio selection with
regime switching: A continuous-time model,
{\sl SIAM J. Control Optim.}, {\bf 42} (2003),
1466--1482.

\bibitem{ZY09}
Zhu, C. and Yin, G.
On competitive Lotka-Volterra model in random environments,
{\sl J. Math. Anal. Appl.},
{\bf 357}
(2009), 154--170.
}\end{thebibliography}
\end{document}